\documentclass[11pt]{article}

\usepackage{latexsym}
\usepackage{amsfonts}
\usepackage{amsmath}

\usepackage[dvips]{color}

\newtheorem{thrm}{Theorem}[section]
\newtheorem{lemma}[thrm]{Lemma}
\newtheorem{prop}[thrm]{Proposition}
\newtheorem{defn}[thrm]{Definition}

\newtheorem{remark}[thrm]{Remark}

\numberwithin{equation}{section}

\usepackage{enumerate}

\usepackage{amssymb}
\usepackage{mathtools}
\mathtoolsset{showonlyrefs}
\usepackage{mathrsfs}
\usepackage{comment}
\usepackage{hyperref}
\usepackage{xcolor}

\def\E{\mathbb{E} }
\def\P{\mathbb{P} }
\def\Q{\mathbb{Q} }
\def\R{\mathbb{R} }
\def\N{\mathbb{N} }

\makeatletter
\begin{document}
\allowdisplaybreaks

\title{\Large \bf{
Local properties for $1$-dimensional critical branching L\'{e}vy process}
	\footnote{The research of this project is supported
     by the National Key R\&D Program of China (No. 2020YFA0712900).}}
\author{ \bf  Haojie Hou \footnote{The research of this author is supported by the China Postdoctoral Science Foundation (No. 2024M764112).
} \hspace{1mm}\hspace{1mm}
Yan-Xia Ren\footnote{The research of this author is supported by NSFC (Grant Nos. 12071011 and 12231002) and LMEQF.
 } \hspace{1mm}\hspace{1mm} and \hspace{1mm}\hspace{1mm}
Renming Song\thanks{Research supported in part by a grant from the Simons
Foundation
(\#960480, Renming Song).}
\hspace{1mm} }
\date{}
\maketitle

\begin{abstract}
Consider a one dimensional critical branching L\'{e}vy process $((Z_t)_{t\geq 0}, \mathbb {P}_x)$. Assume that the offspring distribution either has finite second moment or belongs to the domain of attraction to some $\alpha$-stable distribution
with $\alpha\in (1, 2)$,
and that the underlying  L\'{e}vy process
$(\xi_t)_{t\geq 0}$ is
non-lattice and has finite $2+\delta^*$ moment for some $\delta^*>0$.
We first prove that
$$t^{\frac{1}{\alpha-1}}\left(1- \E_{\sqrt{t}y}\left(\exp\left\{-\frac{1}{t^{\frac{1}{\alpha-1}-\frac{1}{2}}}\int h(x) Z_t(\mathrm{d}x) -\frac{1}{t^{\frac{1}{\alpha-1}}} \int g\left(\frac{x}{\sqrt{t}}\right)Z_t(\mathrm{d}x)\right\}\right)\right)$$
 converges as $t\to\infty$
 for any non-negative bounded  Lipschtitz  function $g$ and any non-negative directly Riemann integrable function $h$ of compact support.
 Then for any $y\in \R$ and
bounded Borel set $A$
 of positive Lebesgue measure with its boundary
having
zero Lebesgue measure,
under a higher moment condition on $\xi$,
we find the decay rate of the probability $\P_{\sqrt{t}y}(Z_t(A)>0)$. As an application, we prove  some convergence results for $Z_t$ under the conditional law $\P_{\sqrt{t}y}(\cdot| Z_t(A)>0).$
 \end{abstract}

\medskip


\medskip


\section{Introduction and main results}

\subsection{Background introduction and motivation}

A branching random walk is a discrete-time Markov process which can be described as follows.
At time $0$ there is a particle at $\mathbf{x}\in \R^d$. At time $1$, this particle dies and gives birth to $N$ offspring with
$\P(N=k ) =p_k$ for $k\in \N:=\{0,1,...\}$, and the relative positions of the offspring to the parent are given by iid copies of
a random variable $\mathbf{X}$. These  offspring form generation 1. Given the information at time 1, at  time 2, individuals
of generation 1 independently repeat their parent's behavior. The procedure goes on. Let $Z_n$ be the counting measure of the individuals of generation $n$. $(Z_n)_{n\ge 0}$ is called a branching random walk starting from an initial individual located at $\mathbf{x}$.  We will use $\P_{\mathbf{x}}$ to denote the law of the branching random walk and $\E_{\mathbf{x}}$ to denote the corresponding expectation.

Let $C^+(\R)$ be the family of non-negative continuous functions on $\R$ and $C_c^+(\R)$ be the subfamily of functions in  $C^+(\R)$ with compact support.
For any $f\in C_c^+(\R)$, we write  $\ell(f):= \int f(x)\mathrm{d}x$.
Assume that the branching random walk is critical, that is,
$\sum_{k=0}^\infty kp_k =1$ and $p_1<1$.
It is well-known that this process will become extinct with probability $1$.
For any $\mathbf{x}\in \R^d$, we use $\| \mathbf{x}\|$ to denote the Euclidean norm.
When $d\ge 3$,  under the assumption
\begin{align}\label{Second-moment}
\sum_{k=0}^\infty k^2 p_k <\infty
\end{align}
and that  either $\mathbf{X}$ is a standard $\R^d$-valued Gaussian random variable or $\mathbf{X}$ is a
bounded symmetric $\mathbb{Z}^d$-valued random variable,
Rapenne  \cite[Lemma 2.10]{Rapenne} proved
that  for any closed ball $\mathbf{A}\subset \R^d$, there exists a constant $I_{\mathbf{A}}$ such that for all
 $\mathbf{x}\in \mathbb Z^d$,
as $n\to\infty$,
\begin{align}\label{Survival-probability}
	 \lim_{n\to\infty} n^{d/2}\P_{[\sqrt{n}]\mathbf{x}} \left(Z_n(\mathbf{A}) > 0\right)  = \frac{I_{\mathbf{A}}}{\sqrt{\mbox{det}(\Sigma)}}\exp\left\{-\frac{1}{2}\mathbf{x}^{T}\Sigma^{-1} \mathbf{x}\right\},
\end{align}
where $\Sigma= (\Sigma_{i,j})_{1\leq i,j\leq d}$ is the covariance
matrix
 of $\mathbf{X}$, i.e., $\Sigma_{i,j}= \mbox{Cov}(X_i, X_j)$ for all $1\leq i,j \leq d$.
Moreover, under the same assumption,
 Rapenne \cite[Proposition 2.13]{Rapenne}
proved that for any
$\mathbf{a}\in \mathbb{Z}^d$
and  closed ball $\mathbf{A}\subset \R^d$, there exists a random point process $(\mathcal{N}_{\mathbf{A}}, \P)$ supported in $\mathbf{A}$ and independent of $\mathbf{a}$ such that
\begin{align}\label{Yaglom-type}
	\P_{[\sqrt{n}]\mathbf{a}}
	\left(Z_n \in \cdot \big| Z_n(\mathbf{A})>0 \right) \stackrel{\mathrm{d}}{\Longrightarrow} \P(\mathcal{N}_{\mathbf{A}}\in \cdot ).
\end{align}
For critical branching Brownian motions and critical super-Brownian motions in dimension $d\geq 3$, results similar to
\eqref{Survival-probability} and \eqref{Yaglom-type} are consequences of \cite[(2.8)]{BCG1997}.
More precisely, taking $f= \theta 1_A$ in \cite[(2.8)]{BCG1997} and then letting $\theta \to \infty$, we get \eqref{Survival-probability}; taking a general  $f\in C_c^+(\R)$ and combining it with \eqref{Survival-probability}, we get \eqref{Yaglom-type}.

When $d=2$, things are quite different. When \eqref{Second-moment} holds and $\mathbf{X}$ is a
$\mathbb{Z}^2$-valued
random variable such that $\P(\| \mathbf{X}\| \leq 1 )=1$,
Lalley and Zheng \cite[Propositions 31 and 33]{LZ11}  proved
that
 for any $\mathbf{x}\in \mathbb{Z}^2$, there exists $C(\mathbf{x})>0$ such that for all $n\geq 2$,
\begin{align}\label{Survival-d=2}
	\frac{1}{C(\mathbf{x})}\leq n(\log n)\P_{[\sqrt{n}]\mathbf{x}}(Z_n(\{\mathbf{0}\})>0) \leq C(\mathbf{x}).
\end{align}
Recently Bai et al \cite{Chen}\footnote{In the journal version, we cited the unpublished preprint, which did not include Bai as a coauthor. We have now found the recent arXiv version of the preprint (see arXiv: 2609.15538 in \cite{Chen}) after checking with the authors, have updated the citation here.} refined the result of \eqref{Survival-d=2}. They proved that if $\sum_{k=0}^\infty k^2 p_k<\infty$, then  for any  $\mathbf{x} \in  \mathbb{Z}^2$,
\begin{align}\label{Refined-Srivival-d=2}
	\lim_{n\to\infty} n(\log n)\P_{[\sqrt{n}]\mathbf{x}}(Z_n(\{\mathbf{0}\})>0) =\frac{4}{\sigma^2}e^{-\frac{5}{4}|\mathbf{x}|^2},
\end{align}
where $\sigma^2:= \sum_{k=0}^\infty k^2 p_k -1$.
Comparing \eqref{Survival-probability} and \eqref{Refined-Srivival-d=2},
we see that there is an extra factor $\log n$ in $d=2$.  In the case of
critical continuous-time binary branching random walk $(Z_t)_{t\ge 0}$ with branching rate $2$, under a second moment condition on the random walk,
Durrett \cite[(8.12)]{Durrett1979} proved that for any bounded open set $\mathbf{A}\subset \R^2$ with
$\ell (\partial \mathbf{A} )=0$
and any $\theta>0$,
\begin{align}\label{eq:d=2}
	\lim_{t\to\infty} t(\log t) \left( 1- \E \left(\exp\left\{-\frac{\theta }{\log t}\frac{Z_t(\mathbf{A})}{\ell (\mathbf{A})}\right\}\right) \right) =\frac{4\theta }{\theta +8\pi },
	\end{align}
 where $\ell$ is the Lebesgue measure.
As a consequence of \eqref{eq:d=2}, see  \cite[(8.15)]{Durrett1979},  for any bounded open set $\mathbf{A}\subset \R^2$ with
$\ell (\partial \mathbf{A} )=0$ and  any $h>0$, it holds that
\begin{align}\label{d=2}
    \lim_{t\to\infty} t (\log t )\P\left(\frac{1 }{\log t} \frac{Z_t(\mathbf{A})}{\ell(\mathbf{A})} >h\right) = 4e^{-8\pi h}.
\end{align}

It the case $d=1$, it is well-known (for example, see
the paragraph below  \cite[Theorem 3]{LS15}, or \cite{HL2024}) that, if the assumption \eqref{Second-moment}  holds
and $\mathbf{E}(X)=0, \mathbf{E}(X^2)=1$, there exists a measure-valued random variable $(Y, \P)$ such that as $n\to\infty$,
\begin{align}\label{Conditioned-law-BRW}
	\P\Big(Z_1^{(n)}\in \cdot \ \big| Z_n (\R)>0\Big)
	\stackrel{\mathrm{d}}{\Longrightarrow} \P(Y\in \cdot),
\end{align}
where $Z_1^{(n)}$ is the random measure such that $\int f(x) Z_1^{(n)}(\mathrm{d}x):= \frac{1}{n}\int f\big(\frac{x}{\sqrt{n}}\big) Z_n(\mathrm{d} x)$ for all   bounded non-negative functions $f$.
The random measure $Y$ is related to
super-Brownian motion, which will be introduced later.
It is easy to see from \eqref{Conditioned-law-BRW} that for any bounded non-negative continuous function $f$ on $\R$,
\begin{align}\label{Eq: Conditioned-law-BRW}
		&\lim_{n\to\infty}
	 \frac{1}{\P(Z_n(\R)>0)} \left(1- \E\left(\exp\left\{-\frac{1}{n}\int f\left(\frac{x}{\sqrt{n}}\right)Z_n(\mathrm{d}x) \right\}\right)\right) \nonumber\\
	&= 1- \E\left(\exp\left\{-\int f(x)Y(\mathrm{d}x)\right\}\right).
\end{align}
We will prove
a result more general than \eqref{Eq: Conditioned-law-BRW}
in Theorem \ref{thm1} below
in the continuous time setting.
To the best of our knowledge, there are no
$d=1$ counterparts yet to the high dimensional results
 \eqref{Survival-probability} and \eqref{Yaglom-type}.
In this paper, we will prove the counterparts  of \eqref{Survival-probability} and \eqref{Yaglom-type} for
1-dimensional critical branching L\'evy processes,
see  Theorems \ref{thm2} and \ref{thm4}  below.

A branching L\'evy process is a continuous counterpart of branching random walk and it can be described as follows. At time 0, there is an
 individual at $\mathbf{x}\in \R^d$ and it moves according to a L\'{e}vy process
 $(\xi_t, \mathbf{P}_\mathbf{x})$.
After an exponential time with parameter $\beta>0$, this individual dies and gives birth to $k$ offspring with probability $p_k$, $k=0, 1, \dots$ located at the parent's death place. The offspring then independently  repeat the parent's behavior.
This procedure goes on.
 Let $Z_t$ be the point process formed by the individual alive at time $t$. The process $(Z_t)_{t\ge 0}$  is called a branching L\'{e}vy process.

 A more general class of branching L\'{e}vy processes on the real line has been introduced
   by Bertoin and Mallein \cite{BM}. Recently, quite a few papers, see, for instance, \cite{BM2018, KV2024, MS2023, SW2019}, have been devoted to this general model.
  In this general model
 the branching rate $\beta$ is allowed to be infinite and that the branching is allowed to be non-local, i.e., the birthplace of the children is not necessarily the deathplace of their parent.
  The branching L\'{e}vy process dealt with in this paper is a special case of the branching L\'{e}vy process introduced in \cite{BM} with local branching and finite branching rate.

We will use $\P_x$ to denote the law of
$(Z_t)_{t\ge 0}$ starting from one individual at $x$
and $\E_x$ to denote the corresponding expectation.
For simplicity, we write $\P:=\P_0$ and $\E:=\E_0$.
We will assume that
$(Z_t)_{t\ge 0}$
is critical:
\[
\sum_{k=0}^\infty kp_k =1\ \mbox{and}\ p_1< 1.
\]

 The main purpose of this paper is to study the asymptotic  behavior of
 $(Z_t)_{t\ge 0}$
 when $d=1$
 under some conditions.
We will assume that
\begin{itemize}
	\item [{\bf(H1)}]  The offspring distribution $\{p_k: k\geq 0\}$ belongs to the domain of attraction of an $\alpha$-stable, $\alpha\in (1,2]$, distribution.
	More precisely, either there exist $\alpha\in (1,2)$ and $\kappa(\alpha)\in(0,\infty)$ such that
	\[
	\lim_{n\to\infty} n^\alpha \sum_{k=n}^\infty p_k = \kappa(\alpha),
	\]
	or that (corresponding $\alpha=2$)
	\[
	\sum_{k=0}^\infty k^2 p_k<\infty.
	\]
\end{itemize}
Under the assumption {\bf(H1)},  it is known (see, for example, \cite{Kolmogorov38, Slack1968, Zolotarev1957}) that
\begin{align}\label{Survival-prob-zeta}
	\lim_{t\to\infty} t^{\frac{1}{\alpha-1}}\P(Z_t(\R)>0 ) =
		\left\{\begin{array}{ll}
		\frac{\alpha-1}{\beta \kappa(\alpha)\Gamma(2-\alpha)},\quad &\mbox{when}\  \alpha\in (1,2),\\
		\frac{2}{\beta \left(\sum_{k=1}^\infty k(k-1)p_k\right)},\quad &\alpha=2,
	\end{array}\right.
\end{align}
For the L\'evy process $(\xi_t)_{t\ge 0}$, we will assume
that
\begin{itemize}
	\item [{\bf(H2)}]
	\begin{align}
		\mathbf{E}_0 (\xi_1)=0,\quad 		\mathbf{E}_0(\xi_1^2)=1;
	\end{align}
\end{itemize}
\begin{itemize}
	\item [\bf(H3)]
	 the law of $\xi_1$ under $\mathbf{P}_0$ is non-lattice;
\end{itemize}
and that
\begin{itemize}
	\item [\bf(H4)]
     there exists
	$\delta^*>0$ such that $\mathbf{E}_0(|\xi_1|^{2+\delta^*})<\infty$.
\end{itemize}

The hypothesis {\bf(H3)}  and {\bf(H4)}
will only be used to prove Lemma \ref{LLT-functional-xi} below.
For some results, we will also need the following  stronger moment condition on the L\'evy process:
\begin{itemize}
	 \item [{\bf(H4')}]
	For the $\alpha \in (1,2]$ in {\bf(H1)}, it holds that
	\begin{align}
		\mathbf{E}_0\left(|\xi_1|^{r_0} \right)<\infty\quad\mbox{for some} \quad r_0 >\frac{2 \alpha}{\alpha-1}.
	\end{align}
\end{itemize}
When $\alpha=2$, the assumption is the same as that in \cite{LS15}.
 For any $t>0$, define $M_t$ to be the maximal position of all the
particles at time $t$.
We
 also define
\[
M:= \sup_{t>0} M_t.
\]
Under {\bf(H1)}, {\bf(H2)} and the weaker moment condition $\mathbf{E}_0((\xi_1 \vee 0)^{r_0} )<\infty$  than
{\bf(H4')},
it was proved in \cite{HJRS} (although \cite{HJRS} did not deal with  the case $\alpha=2$,  the proof is actually the same as the case $\alpha\in (1,2)$, see the argument below \cite[Theorem 1.1]{HJRS}) that there exists $\theta (\alpha)\in (0,\infty)$ such that
\begin{align}\label{Tail-of-M}
	\lim_{x\to\infty} x^{\frac{2}{\alpha-1}}\P(M\geq x) = \theta (\alpha).
\end{align}
The assumption {\bf(H4')}
 is only used in the proof of Lemma \ref{Technical-Lemma} to control
 the  overshoot of the
underlying L\'evy process.

\subsection{Critical super-Brownian motion}

In this subsection, we give a brief introduction to super-Brownian motion.
Let $\mathcal{M}_F(\R)$ be the families of finite Borel measures on $\R$. We will use $\mathbf{0}$ to denote  the null measure on $\R$. Let $B_b^+(\R)$ be the space of  non-negative bounded Borel functions on $\R$. For any $f\in B_b^+(\R)$ and $\mu\in \mathcal{M}_F(\R)$, we use $\mu (f)$ to denote the integral of $f$ with respect to $\mu$. For any $\alpha\in (1, 2]$, the function
\begin{align}\label{Stable-Branching-mechanism}
	\varphi(\lambda) := \mathcal{C}(\alpha) \lambda^\alpha:=	\left\{\begin{array}{ll}
		\frac{\beta \kappa(\alpha)\Gamma(2-\alpha)}{\alpha-1}  \lambda^\alpha,\quad &\mbox{when}\  \alpha\in (1,2),\\
		\frac{\beta }{2}\left(\sum_{k=1}^\infty k(k-1)p_k\right) \lambda^2,\quad &\alpha=2,
	\end{array}\right.
\end{align}
where $\kappa(\alpha)$ is the constant
given in {\bf(H1)} and $\Gamma(z):=\int_0^\infty t^{z-1}e^{-t}\mathrm{d}t$ is the Gamma function,  is a branching mechanism. Since $\varphi'(0)=0$, $\varphi$ is a critical branching mechanism.
Let $(B_t, \mathbf{P}_x)$ be a standard Brownian motion.

The critical super-Brownian motion $X=\{(X_t)_{t\geq 0}; \mathbb P_\mu\}$ that we will use in this paper
is an $\mathcal{M}_F(\R)$-valued Markov process such that for any $f\in B_b^+(\R)$,
$$
-\log  \E_{\mu} \left(\exp\left\{ -X_t(f)\right\}\right)=\mu \left(v_f^X(t, \cdot)\right),
$$
where $(t, x)\mapsto v_f^X(t, x)$ is the unique locally bounded non-negative solution to
\begin{align}\label{Evolution-cumulant-semigroup}
	v^X_f(t, x)&=\mathbf{E}_x\left(f(B_t) \right)- \mathbf{E}_y\Big(\int_0^t \varphi\left( v_f^X(t-s,B_s)\right)\mathrm{d}s\Big).
\end{align}
Since $1< \frac{2}{\alpha-1}$, by \cite[Theorem 1.2]{FMW2010} and \cite{KS1988}, for any $\mu\in \mathcal{M}_F(\R)$, $\P_\mu$-almost surely, the random measure $X_t$ is absolutely continuous with respect to the Lebesgue measure and the density function
\[
Y_t(x):= \frac{X_t(\mathrm{d}x)}{\mathrm{d} x}
\]
has a version which is continuous in $x$. We will always use  $Y_t$ to denote this version.

For the probabilistic representation of the weak convergence limit via super Brownian motion in
Theorem  \ref{thm4}
below, we will also need the $\N$-measures of super Brownian motion.

Without loss of generality, we assume that $X$ is the coordinate process on
\[
\mathbb D:=\{ w= (w_t)_{t\geq 0}: w \text{ is an $\mathcal{M}_F(\R)$-valued	c\`{a}dl\`{a}g	function} \}.
\]
We assume that $(\mathcal{F}_\infty, (\mathcal{F}_t)_{t\ge 0})$ is the natural filtration on $\mathbb D$, completed as usual with the $\mathcal{F}_\infty$-measurable and $\mathbb P_\mu$-negligible sets for all $\mu\in\mathcal{M}_F(\R)$. Let $\mathbb W^+_0$ be the family of
$\mathcal{M}_F(\R)$-valued c\`{a}dl\`{a}g functions on $(0, \infty)$ with $\mathbf{0}$ as a trap and with
$\lim_{t\downarrow0}w_t= \mathbf{0}$.

Since the super Brownian motion $X_t$ is critical and that $\int_1^\infty \frac{1}{\varphi(\lambda)}\mathrm{d}\lambda<\infty$, we see that $\P_{\delta_y} (X_t=\mathbf{0}) >0$ for all $t>0$ and  $y\in \R$,  which implies that there exists a unique family of $\sigma$-finite measures $\{\mathbb N_y; y\in \R \}$ on $\mathbb W^+_0$ such that for any $\mu\in \mathcal {M}_F(\R)$,
if ${\mathcal N}(\mathrm{d}w)$ is a Poisson random measure on $\mathbb W^+_0$ with intensity measure
$$
\mathbb N_\mu(\mathrm{d}w):=\int_{\R} \mathbb N_y(\mathrm{d}w)\mu(\mathrm{d}y),
$$
then the process defined by
$$
\widehat X_0:=\mu, \quad \widehat X_t:=\int_{\mathbb W^+_0 }w_t{\mathcal N}(\mathrm{d}w), \quad t>0,
$$
is a realization of the superprocess $X=\{(X_t)_{t\geq 0}; \mathbb P_\mu\}$.
Furthermore, for any $t>0$, $y\in \R$ and  $f\in B^+_b(\R)$,
\begin{align}\label{N-measure-equation}
	\mathbb N_y\left(1- \exp\left\{-
	w_t(f)
	\right\}\right)= -\log \mathbb{E}_{\delta_y}\left(\exp\left\{-
	X_t(f)\right\}\right),
\end{align}
see \cite{DyKu} or \cite[Theorems 8.27 and 8.28]{LZ}.
The next useful result says that for any given $t>0$ and $y\in\R$, $w_t$ has an $\N_y$-a.e. continuous
density.

Define
\[
\mathcal{A}:= \left\{\mu\in \mathcal{M}_F(\R):  \frac{\mathrm{d}\mu}{\mathrm{d}x} \in C^+(\R) \right\}.
\]

\begin{lemma}\label{Lemma:Density-under-N-measure}
	For any $t>0$ and $y\in \R$, it holds that
	\[
	\N_y \left( w_t\notin \mathcal{A}\right)=0.
	\]
\end{lemma}
The proof is postponed to Section \ref{Section:Lemma:Density-under-N}.  We still use $\{Y_t(x), x\in\R\}$ to denote
the density of $w_t$.

\subsection{Main results}

We will sometimes use $\ell(A)$ to denote the Lebesgue measure of a Borel set $A \subset\R$.
 Let ${\rm DRI}^+(\R)$ (${\rm DRI}^+_c(\R)$) be the family of
 non-negative directly Riemann integrable
 functions (of compact support).
We say that a bounded Borel set $A$ is
directly Riemann integrable  if the indicator $1_A$ is a directly Riemann integrable function.
It is well known that (i) any directly Riemann integrable function is bounded; (ii) a non-negative Borel function of compact support is directly Riemann integrable if and only if it is
Riemann integrable and (iii) a bounded Borel set $A$ is directly Riemann integrable if and only if $\ell(\partial A)=0$.
For the definition of  directly Riemann integrable function,
see the beginning of Subsection \ref{ss:2.2}.
Let $B_{Lip}^+(\R)$ be the family of bounded non-negative Lipschitz continuous functions in $\R$. For any $g\in B_{Lip}^+(\R)$, we use $\mbox{Lip}(g)$ to denote its Lipschitz constant.

\begin{thrm}\label{thm1}
		 Assume  {\bf(H1)},  {\bf(H2)},  {\bf(H3)} and {\bf(H4)} hold.
		Then for any
		$y\in \R$, $g\in B_{Lip}^+(\R)$ and
				$h\in {\rm DRI}^+_c(\R)$,
		it holds that
	\begin{align}\label{Limit}
		& \lim_{t\to\infty} t^{\frac{1}{\alpha-1}}\left(1- \E_{\sqrt{t}y}\left(\exp\left\{-\frac{1}{t^{\frac{1}{\alpha-1}-\frac{1}{2}}}\int h(x) Z_t(\mathrm{d}x) -\frac{1}{t^{\frac{1}{\alpha-1}}} \int g\left(\frac{x}{\sqrt{t}}\right)Z_t(\mathrm{d}x)\right\}\right)\right)\nonumber\\
		&= -\log \E_{\delta_y}\left(\exp\left\{- \ell(h) Y_1(0) -X_1(g)\right\}\right).
	\end{align}
\end{thrm}

\begin{remark}\label{Remark}
	In the special case $\alpha=2$,
		taking $h=0$ in Theorem \ref{thm1}, we get that
	\begin{align}\label{Eq:Remark1}
			& \lim_{t\to\infty} t\left(1- \E_{\sqrt{t}y}\left(\exp\left\{ -\frac{1}{t} \int g\left(\frac{x}{\sqrt{t}}\right)Z_t(\mathrm{d}x)\right\}\right)\right)\nonumber\\
		&= -\log \E_{\delta_y}\left(\exp\left\{- X_1(g)\right\}\right)= \N_y\left(1-\exp\left\{-w_1(g)  \right\}\right),
	\end{align}
where in the last equality we used \eqref{N-measure-equation}.
	Combining \eqref{Survival-prob-zeta} and \eqref{Eq:Remark1}, we get that
	\begin{align}\label{Eq:Remark2}
			&\lim_{t\to\infty}  \frac{1}{\P(Z_t(\R)>0)} \left(1- 
			\E_{\sqrt{t}y}
			\left(\exp\left\{-\frac{1}{t}\int g\left(\frac{x}{\sqrt{t}}\right)Z_t(\mathrm{d}x) \right\}\right)\right) \nonumber\\
			&=
						\mathcal{C}(2)
			\N_y\left(1-\exp\left\{-w_1(g)  \right\}\right).
						\end{align}
It follows from \eqref{N-measure-equation} and \cite[(1.11)]{HRS}
that $ \N_y(w(1)>0)=\mathcal{C}(2)^{-1}$.
Therefore, by \eqref{Eq:Remark2}, we conclude that
	\begin{align}\label{Eq:Remark3}
			&\lim_{t\to\infty}  \frac{1}{\P(Z_t(\R)>0)} \left(1- 
			\E_{\sqrt{t}y}
			\left(\exp\left\{-\frac{1}{t}\int g\left(\frac{x}{\sqrt{t}}\right)Z_t(\mathrm{d}x) \right\}\right)\right) \nonumber\\
			&= \N_y\left(1-\exp\left\{-w_1(g)  \right\}|w(1)>0\right).
	\end{align}
	Combining \eqref{Eq: Conditioned-law-BRW} and \eqref{Eq:Remark3}, we immediately get that
	$(Y, \P)\stackrel{\mathrm{d}}{=} (w_1, \N_0(\cdot| w_1(1)>0))$.

In the special case $\alpha=2$, taking $g=0$ and
$h(x)=\frac{\theta 1_A(x)}{\ell(A)}$  with $\theta >0$,
 $A$ being a bounded Borel set
with $\ell(A)>0$ and $\ell(\partial A)=0$ in Theorem \ref{thm1}, we get
\begin{align}
\lim_{t\to\infty}t\left(1- \E\left(\exp\left\{-
\frac{\theta }{\sqrt{t}}
\frac{Z_t(A)}{\ell(A)} \right\}\right)\right)=-\log \E_{\delta_y}\left(\exp\left\{-
\theta Y_1(0)\right\}\right).
\end{align}
Comparing the result above with  \eqref{eq:d=2} and \eqref{d=2} for $d=2$,
we see the differences between the cases $d=2$ and $d=1$. In the case $d=2$,  there is an extra factor $\log t$ in the decay, and also one
needs
 to normalize with  $\log t$ instead of $\sqrt{t}$. In addition, the limit in $d=2$ is related to the Laplace transform of an exponential random variable while
  in  the case $d=1$,
 the limit is  related to super-Brownian motion.
\end{remark}

\begin{thrm}\label{thm2}
		Assume  {\bf(H1)},  {\bf(H2)},  {\bf(H3)} and {\bf(H4')} hold.
			Then for any $y\in \R$ and any  bounded Borel set $A$ with  $\ell(A)>0$ and $\ell(\partial A)=0$,  it holds that
	\begin{align}
		 \lim_{t\to\infty} t^{\frac{1}{\alpha-1}}\P_{\sqrt{t}y}\left(Z_t(A)>0 \right)= -\log \P_{\delta_y}\left( Y_1(0)=0\right).
	\end{align}
\end{thrm}
  \begin{remark}
  	  	(i) Taking $h=\theta 1_A$ and $g=0$ in Theorem \ref{thm1} first and then letting $\theta\uparrow \infty$,
 formally one expect
the limit in Theorem \ref{thm2} to be $-\log \P_{\delta_y}\left( \ell(A)Y_1(0)=0\right).$
Under the assumption  $\ell(A)>0$, the limit is equal to
    $-\log \P_{\delta_y}\left(Y_1(0)=0\right)$,
which implies that  the right-hand side of the limit in Theorem \ref{thm2} is independent of $A$.

  	(ii)
   	 When $\alpha=2$, Theorem \ref{thm2}  is the
	 	 1-d counterpart to the high dimensional result  \eqref{Survival-probability}
	 	 and \eqref{Refined-Srivival-d=2}.
	     We see that
		 branching plays a more important role in dimension 1 while spatial motion dominates in dimension
		 $d\geq 2$.
	 	 In dimension 1, the limit is related to the density of super-Brownian motion, while in dimension  $d\geq 3$, 	the limit in \eqref{Survival-probability} is only related to the local limit  of a random walk (see \cite[Proposition 2.1]{Rapenne}) and that branching only appears in the constant $I_A$ (see the end of the proof of \cite[Lemma 2.10]{Rapenne} on page 14). In dimension 2,
		 both the branching and the spatial motion
		 		 affect
		 the limit in \eqref{Refined-Srivival-d=2} in the sense that the limit requires at least second moment due to the appearance of $\sigma^2$ and that the exponential term $e^{-\frac{5}{4}|x|^2}$ is related to the local limit theorem for the random walk.
 \end{remark}

For any $t>0$, we define a measure $Z_1^{(t)}$ by
\[
\int f(y)Z_1^{(t)}(\mathrm{d}y):= \frac{1}{t^{\frac{1}{\alpha-1}}}\int f\left(\frac{y}{\sqrt{t}}\right) Z_t(\mathrm{d}y).
\]
The next result is an application of Theorems \ref{thm1} and  \ref{thm2}.

\begin{thrm}\label{thm4}
		Assume that {\bf(H1)},  {\bf(H2)},  {\bf(H3)} and {\bf(H4')} hold. Suppose that
		$y\in\R$ and $A$ is a bounded Borel set with $\ell(A)>0$ and $\ell(\partial A)=0$.
	
	(i) As $t\to\infty$, we have
    \begin{align}
    	\left(\frac{1}{t^{\frac{1}{\alpha-1}-\frac{1}{2}}}Z_t,\P_{\sqrt{t}y}(\cdot| Z_t(A)>0)\right) \stackrel{\mathrm{d}}{\Longrightarrow} \left( Y_1(0)   \ell  , \N_y(\cdot| Y_1(0)>0) \right)
    \end{align}
     in the sense of vague  topology.

     (ii) As $t\to\infty$, it holds that
     \begin{align}
     	 \left(Z_1^{(t)}, \P_{\sqrt{t}y}(\cdot| Z_t(A)>0) \right) \stackrel{\mathrm{d}}{\Longrightarrow} \left(w_1, \N_y(\cdot|Y_1(0)>0)\right)
     \end{align}
     in the sense of weak topology.
\end{thrm}
\begin{remark}
		In the special case $\alpha=2$, Theorem \ref{thm4} (i) is the 1-d counterpart to the high dimensional result  \eqref{Yaglom-type}.  There are some differences between the 1-d case and the high dimensional case. First there is an extra factor $\sqrt{t}$ in the 1-d case while no normalization in the high dimensional case $d\ge 3$. Also in the 1-d case, the limit is
an absolutely continuous random measure (with respect to the Lebesgue measure) with density $Y_1(0)$ (the density $Y_1(x)$ of super-Brownian motion $X_1$ evaluated at 0),
while in the high dimensional case $d\ge 3$, the limit $\mathcal{N}_A$ is
a random point measure supported on $A$.

	Theorem \ref{thm4} (ii) should be compared with \eqref{Conditioned-law-BRW}. \eqref{Conditioned-law-BRW} is about the asymptotic of  $Z_t$ conditioned on  global survival $Z_t(\R)>0$, while Theorem \ref{thm4} (ii)
is about the asymptotic of  $Z_t$ conditioned on  local $Z_t(A)>0$.
As we mentioned in Remark \ref{Remark}, in the special case  $\alpha=2$,
the limit $(Y, \P)$ in \eqref{Conditioned-law-BRW} is  equal in law to $(w_1, \N_0(\cdot| w_1(1)>0))$ which is different from the limit $ \left(w_1, \N_y(\cdot|Y_1(0)>0)\right)$  in Theorem \ref{thm4} (ii).
	
	Theorem \ref{thm4} (i) describes the local behavior of the counting measure $Z_t$, while Theorem \ref{thm4} (ii)
	is about the global behavior of $Z_t$.
\end{remark}

After our paper was accepted, we found that Lee \cite{Lee} has already studied local properties of critical binary branching Brownian
motion in all dimensions $d\geq 1$. In the present paper, our spatial motion is a genaral L\'evy process and our branching mechanism is
only assumed to be in domain of attaction of a stable distribution. Moreover, our proof strategies are different from \cite{Lee}.

We end this section with a brief description of the organization of this paper.
In Section \ref{S:2}, we give some elementary estimates involving the standard normal density and about the underlying L\'evy process. We also derive an integral equation for the Laplace transform of $Z_t$ and prove the existence and  uniqueness of solution
for the problem \eqref{PDE} below.
In Section  \ref{S:3}, we give the proofs of  Theorems \ref{thm1}, \ref{thm2} and \ref{thm4}.
In Section \ref{Section:Lemma:Density-under-N}, we give the proof of Lemma \ref{Lemma:Density-under-N-measure}.

For two functions $f(x)$ and $g(x)$ with $x\in E$, we use $f\lesssim g, x\in E,$ to denote that there exists a constant $C$ independent of $x$ such that $f(x)\leq C g(x), x\in E$.

\section{Preliminaries}\label{S:2}

\subsection{Some estimates involving  the standard normal density}

Throughout this paper,
$\phi(x):= \frac{1}{\sqrt{2\pi}}e^{-\frac{x^2}{2}}$ is
the  standard normal density.

\begin{lemma}\label{lemma2}
	\begin{enumerate}
		\item[(i)] 	For any $\Delta>0$,
		\begin{align}
			\sup_{y\in \R} \left| \phi(y)- \phi(y+\Delta)\right|\leq (\Delta \wedge \sqrt{\Delta}).
		\end{align}
		\item [(ii)] For any $0<r<s$ with $s-r\in (0,1)$, it holds that
		\begin{align}
			\sup_{y\in\R} \left| \frac{1}{\sqrt{r}}\phi\left(\frac{y}{\sqrt{r}}\right)-\frac{1}{\sqrt{s}}\phi\left(\frac{y}{\sqrt{s}}\right)\right| \leq \frac{1}{\sqrt{r}}-\frac{1}{\sqrt{s}} + \frac{1}{\sqrt{s}}\left(\sqrt{s-r}+ 1-\exp\left\{-\frac{\sqrt{s-r}}{r}\right\} \right).
		\end{align}
	\end{enumerate}
\end{lemma}
\textbf{Proof: } (i)
It is easy to check that
\begin{align}
	|\phi'(x)|=\frac{1}{\sqrt{2\pi}}|x|e^{-\frac{x^2}{2}}\leq 1.
\end{align}
Therefore, noticing that $\phi(x)\leq \frac{1}{\sqrt{2\pi}}\leq \frac{1}{2}$, we conclude that
\[
\sup_{y\in\R} \left|\phi(y)-\phi(y+\Delta)\right|\leq \Delta \land 1 \leq
(\Delta \wedge \sqrt{\Delta}).
\]

(ii) For any $y\in \R$,
\begin{align}\label{Goal-2}
	& \left|\frac{1}{\sqrt{r}}\phi\left(\frac{y}{\sqrt{r}}\right)- \frac{1}{\sqrt{s}}\phi\left(\frac{y}{\sqrt{s}}\right)\right|\leq  \left|\frac{1}{\sqrt{r}}-\frac{1}{\sqrt{s}}\right| \phi\left(\frac{y}{\sqrt{r}}\right)+ \frac{1}{\sqrt{s}}\left| \phi\left(\frac{y}{\sqrt{r}}\right)- \phi\left(\frac{y}{\sqrt{s}}\right)\right| \nonumber\\
	& \leq \frac{1}{\sqrt{2\pi}}\left|\frac{1}{\sqrt{r}}-\frac{1}{\sqrt{s}}\right| + \frac{1}{\sqrt{2\pi s}}\exp\left\{-\frac{y^2}{2s}\right\}\left(1- \exp\left\{-y^2\left(\frac{1}{2r} -\frac{1}{2s}\right) \right\} \right) \nonumber\\
	& \leq \frac{1}{\sqrt{r}}-\frac{1}{\sqrt{s}}+ \frac{1}{\sqrt{s}} \exp\left\{-\frac{y^2}{2s}\right\}\left(1- \exp\left\{-y^2\left(\frac{1}{2r} -\frac{1}{2s}\right) \right\} \right).
\end{align}
If $y^2 \sqrt{s-r}>2s$, then by the  inequality $ae^{-a}<1$ for all $a>1$, we get that
\begin{align}\label{Proof-Goal-2-1}
	&\exp\left\{-\frac{y^2}{2s}\right\}\left(1- \exp\left\{-y^2\left(\frac{1}{2r} -\frac{1}{2s}\right) \right\} \right)\leq \exp\left\{-\frac{1}{\sqrt{s-r}}\right\} \leq \sqrt{s-r}.
\end{align}
If $y^2 \sqrt{s-r}\leq 2s$, then
\begin{align}\label{Proof-Goal-2-2}
		&\exp\left\{-\frac{y^2}{2s}\right\}\left(1- \exp\left\{-y^2\left(\frac{1}{2r} -\frac{1}{2s}\right) \right\} \right)\leq 1-\exp\left\{ -\frac{s}{\sqrt{s-r}} \cdot \frac{s-r}{sr}\right\}\nonumber\\
		& = 1-\exp\left\{-\frac{\sqrt{s-r}}{r}\right\}.
\end{align}
Combining \eqref{Goal-2}, \eqref{Proof-Goal-2-1} and \eqref{Proof-Goal-2-2}, we conclude that
\begin{align}
	& \left|\frac{1}{\sqrt{r}}\phi\left(\frac{y}{\sqrt{r}}\right)- \frac{1}{\sqrt{s}}\phi\left(\frac{y}{\sqrt{s}}\right)\right|\leq  \frac{1}{\sqrt{r}}-\frac{1}{\sqrt{s}} + \frac{1}{\sqrt{s}}\left(\sqrt{s-r}+ 1-\exp\left\{-\frac{\sqrt{s-r}}{r}\right\} \right),
\end{align}
which implies the desired result.
\hfill$\Box$

The following inequality will be used several times later:
\begin{align}\label{Lemma:inequality}
		\left|1-e^{-(x+y)} -x-y\right|\leq x^2 +y^2, \quad x, y\ge 0.
	\end{align}
The proof of this inequality is elementary and we omit it.

\subsection{Estimates for the L\'{e}vy process}\label{ss:2.2}

We first give a local limit theorem for the underlying L\'evy process
$(\xi_t)_{t\geq 0}$.
Before that, we recall the definition
    of directly Riemann integrable functions. For more details on properties of  directly Riemann integrable functions,
one can refer to \cite[Section XI.1]{Feller} and \cite[Section 2.1]{GX}.

Let $f$ be a non-negative Borel function.  For any $\kappa >0$,  define
\begin{align}
&\overline{f}_\kappa (x):= \sum_{m \in \mathbb{Z}} 1_{[m\kappa, (m+1)\kappa)}(x) \sup_{z\in [m\kappa, (m+1)\kappa)} f(z),\nonumber\\
&\underline{f}_\kappa (x):= \sum_{m \in \mathbb{Z}} 1_{[m\kappa, (m+1)\kappa)} (x)\inf_{z\in [m\kappa, (m+1)\kappa)} f(z).
\end{align}
We say that $f$ is directly Riemann integrable if $\int \overline{f}_\kappa (x) \mathrm{d} x<\infty$ for
some $\kappa>0$ and
\[
\lim_{\kappa\to 0}\int_{\R}\left(\overline{f}_\kappa (x)-\underline{f}_\kappa (x)\right)dx=0.
\]
Recall that we use ${\rm DRI}^+(\R)$ to denote the family of  non-negative directly Riemann integrable functions.
It is easy to see
from the definition that any $h\in  {\rm DRI}^+(\R)$ must be bounded.
For $h\in{\rm DRI}^+(\R)$, we define $\Vert h\Vert_\infty:= \sup_{x\in \R}
|h(x)|$.

\begin{lemma}\label{LLT-functional-xi}
Assume that {\bf(H2)}, {\bf(H3)} and {\bf(H4)} hold.
For any $f\in {\rm DRI}^+_c(\R)$, it holds that
	\begin{align}\label{LLT-for-F}
		\lim_{n\to\infty} \sup_{x\in \R} \left|   \sqrt{n}\mathbf{E}_x\left(f(\xi_n)\right)-	\ell(f)
		\phi\left(\frac{x}{\sqrt{n}}\right)\right| =0,
	\end{align}
where $\phi(x)= \frac{1}{\sqrt{2\pi}}e^{-\frac{x^2}{2}}$ is the standard normal density.
\end{lemma}

\textbf{Proof: }
For any $\kappa, \vartheta>0$, define
\begin{align}
	\overline{f}_{\kappa, \vartheta}(x):= \sup_{|y|\leq \vartheta} \overline{f}_\kappa(x+y)\quad \mbox{and}\quad \underline{f}_{\kappa, \vartheta}(x):= \inf_{|y|\leq \vartheta}
	\underline{f}_\kappa(x+y),
\end{align}
then by \cite[Lemma 2.2]{GX}, it holds that
\begin{align}\label{Uniformly-bounded}
 \lim_{\kappa \to 0} \lim_{\vartheta\to 0} \int \left|\overline{f}_{\kappa, \vartheta}(x) -f(x)\right|\mathrm{d}x = \lim_{\kappa \to 0} \lim_{\vartheta\to 0} \int \left|\underline{f}_{\kappa, \vartheta}(x) -f(x)\right|\mathrm{d}x =0.
\end{align}
Let $\vartheta\in (0, \frac{1}{2})$ be sufficiently small, then by \cite[(2.6) and Theorem 2.7]{GX}, there exist a constant
$K>0$ independent of $\vartheta$
 and $\kappa$ and a constant $C_\vartheta>0$ independent of $\kappa$ such that for any $\kappa>0, x\in \R$,
\begin{align}\label{Upper-bound-LLT}
	\mathbf{E}_0\left(f(x+\xi_n)\right)-
\frac{1+K \vartheta}{\sqrt{n}}
\int\overline{f}_{\kappa, \vartheta} (x+z)
	\phi\left(\frac{z}{\sqrt{n}}\right)\mathrm{d}z \leq \frac{C_\vartheta}{n^{(1+\delta^*)/2}}
	\int  \overline{f}_{\kappa, \vartheta}  (x+z)\mathrm{d}z
\end{align}
and that
\begin{align}\label{Lower-bound-LLT}
		& \mathbf{E}_0\left(f(x+\xi_n)\right)-\frac{1}{\sqrt{n}}\int \left(
		\underline{f}_{\kappa, \vartheta}(x+z)
		-K\vartheta
f(x+z) \right)\phi\left(\frac{z}{\sqrt{n}}\right)\mathrm{d}z\nonumber\\
		&\geq -\frac{C_\vartheta}{
			n^{(1+\delta^*)/2}}
			\int f(x+z)\mathrm{d}z.
\end{align}
Therefore, by \eqref{Upper-bound-LLT}, we see that
\begin{align} \label{Upper-bound-LLT-2}
	& \sqrt{n}\mathbf{E}_x\left(f(\xi_n)\right)- \int   f (x+z)\phi\left(\frac{z}{\sqrt{n}}\right)\mathrm{d}z \nonumber\\
	& \leq \left(
K\vartheta
 +\frac{C_\vartheta}{
		n^{\delta^*/2}}\right)
		\int  \overline{f}_{\kappa, \vartheta} (x+z)\mathrm{d}z+
		\int \left|\overline{f}_{\kappa, \vartheta}(x) -f(x)\right|\mathrm{d}x \nonumber\\
	& \leq \left(K\vartheta +\frac{C_\vartheta}	{n^{\delta^*/2}} +1	\right)	\int \left|\overline{f}_{\kappa, \vartheta}(x) -f(x)\right|\mathrm{d}x + \left(K\vartheta +\frac{C_\vartheta}	{n^{\delta^*/2}} 	\right)	\int f(x)\mathrm{d}x	=: \overline{I}(\vartheta, \kappa, n).
\end{align}
Similarly, according to \eqref{Lower-bound-LLT}, we have
\begin{align}\label{Lower-bound-LLT-2}
	&  \sqrt{n}\mathbf{E}_x\left(f(\xi_n)\right)- \int   f (x+z)\phi\left(\frac{z}{\sqrt{n}}\right)\mathrm{d}z\nonumber\\
	& \geq -\frac{C_\vartheta}	{n^{\delta^*/2}} 	\int f(x+z)\mathrm{d}z-
	K\vartheta
	\int f(x+z)\phi\left(\frac{z}{\sqrt{n}}\right)\mathrm{d}z- \int \left|\underline{f}_{\kappa, \vartheta}(x) -f(x)\right|\mathrm{d}x\nonumber\\
	& \geq -\left(K\vartheta +\frac{C_\vartheta}	{n^{\delta^*/2}} 	\right) \int f(x)\mathrm{d}x	-\int \left|\underline{f}_{\kappa, \vartheta}(x) -f(x)\right|\mathrm{d}x	\geq  -
	 \underline{I}(\vartheta, \kappa, n).
	\end{align}
Therefore, combining \eqref{Upper-bound-LLT-2} and \eqref{Lower-bound-LLT-2}, we conclude that
\begin{align}
	& \limsup_{ n\to\infty}
	\sup_{x\in \R} \left|   \sqrt{n}\mathbf{E}_x\left(f(\xi_n)\right)- \int   f (x+z)\phi\left(\frac{z}{\sqrt{n}}\right)\mathrm{d}z\right| \nonumber\\
		&
    \leq \lim_{ n\to\infty}	\overline{I}(\vartheta, \kappa, n) + \lim_{ n\to\infty}	\underline{I}(\vartheta, \kappa, n)  \nonumber\\
& = \left(K\vartheta +1	\right)	\int \left|\overline{f}_{\kappa, \vartheta}(x) -f(x)\right|\mathrm{d}x +2K\vartheta 	\int f(x)\mathrm{d}x	+\int \left|\underline{f}_{\kappa, \vartheta}(x) -f(x)\right|\mathrm{d}x.
\end{align}
By \eqref{Uniformly-bounded}, letting $\vartheta \to 0$ first and then $\kappa \to 0$ in the above inequality, we get
\begin{align}
		  \lim_{ n\to\infty}
	 \sup_{x\in \R} \left|   \sqrt{n}\mathbf{E}_x\left(f(\xi_n)\right)- \int   f (x+z)\phi\left(\frac{z}{\sqrt{n}}\right)\mathrm{d}z\right| =0.
\end{align}
Thus, to prove the desired result, it remains to show that
\begin{align}\label{step_27}
	 \lim_{n\to\infty}
	 \sup_{x\in \R} \left|  \phi\left(\frac{x}{\sqrt{n}}\right) \int   f (z) \mathrm{d}z - \int   f (x+z)\phi\left(\frac{z}{\sqrt{n}}\right)\mathrm{d}z\right| =0.
\end{align}
Let $E$ be any bounded interval such that $\mbox{supp}(f)\subset E$, then for any $x\in \R$,
\begin{align}
	& \left|  \phi\left(\frac{x}{\sqrt{n}}\right) \int   f (z) \mathrm{d}z - \int   f (x+z)\phi\left(\frac{z}{\sqrt{n}}\right)\mathrm{d}z\right| \nonumber\\
	& =\left|  \phi\left(\frac{x}{\sqrt{n}}\right) \int   f (z) \mathrm{d}z - \int   f (z)\phi\left(\frac{x-z}{\sqrt{n}}\right)\mathrm{d}z\right| \nonumber\\
	& \leq \Vert f\Vert_\infty  \int_E \left|  \phi\left(\frac{x}{\sqrt{n}}\right)  - \phi\left(\frac{x-z}{\sqrt{n}}\right) \right|\mathrm{d}z \leq\ell(E) \Vert f\Vert_\infty  \sup_{z\in E} \left|  \phi\left(\frac{x}{\sqrt{n}}\right)  - \phi\left(\frac{x-z}{\sqrt{n}}\right)\right|.
\end{align}
If $|x|<n^{2/3}$, then by the inequality $e^{-a}-e^{-b}\leq |b-a |$ for any $a, b\geq 0,$
we see that
\begin{align}\label{step_28}
	&\sup_{z\in E} \left|  \phi\left(\frac{x}{\sqrt{n}}\right)  - \phi\left(\frac{x-z}{\sqrt{n}}\right)\right|\leq \frac{1}{\sqrt{2\pi}}\sup_{z\in E} \left| \frac{x^2- (x-z)^2}{n}\right| \nonumber\\
	&\leq  \frac{2}{\sqrt{2\pi}}
	\left( \frac{1}{n}\sup_{z\in E}z^2 + \frac{1}{n^{1/3}} \sup_{z\in E} |z| \right) \stackrel{n\to\infty}{\longrightarrow}0.
\end{align}
On the other hand, if $|x|>n^{2/3}$, then for large $n$, we see that
\begin{align}\label{step_29}
	&\sup_{z\in E} \left|  \phi\left(\frac{x}{\sqrt{n}}\right)  - \phi\left(\frac{x-z}{\sqrt{n}}\right)\right|\leq \phi (n^{1/6})+ \phi\left(\frac{n^{2/3}-\sup_{z\in E}|z|}{\sqrt{n}}\right)\stackrel{n\to\infty}{\longrightarrow}0.
\end{align}
The proof is now complete.
\hfill$\Box$

\begin{remark}
     We mention here that the non-lattice assumption {\bf(H3)} is
      only used to prove Lemma \ref{LLT-functional-xi}.
       If  {\bf(H3)} does not hold, it is possible to get a result similar to  Lemma \ref{LLT-functional-xi}.
     For example,    if $\xi$ is a compound Poisson process
     supported on $\mathbb{Z}$
  with $\mathbf{E}_0(\xi_1)=0, \mathbf{E}_0(\xi_1^2)=1$, $\mathbf{E}_0(|\xi_1|^3)<\infty$ and that the  support of the L\'{e}vy measure  contains $\{n, n+1\}$ for some $n\in \N$,
  then by
  \cite[Theorem 13, p.206]{Petrov1975}, for any $a\in  \mathbb{Z}$, it is easily seen that
     \[
     	\lim_{n\to\infty} \sup_{x\in \mathbb{Z}} \left|   \sqrt{n}\mathbf{P}_x\left(\xi_n=a\right)-	
     \phi\left(\frac{x}{\sqrt{n}}\right)\right|
     =0.
     \]
     Replace
     the Lebesgue measure $\ell$ by the counting measure $\ell_c$ on $\mathbb{Z}$,   and for any bounded function $f$ with compact support, define
     \[
     \ell_c(f):= \sum_{i\in\mathbb{Z}} f(i).
     \]
    Denote by $B_c^+(\mathbb{Z})$ the class of non-negative bounded functions with compact support. In this case, we see that for any $f \in B_c^+(\mathbb{Z})$,
      \[
     \lim_{n\to\infty} \sup_{x\in \mathbb{Z}} \left|   \sqrt{n}
     \mathbf{E}_x
     \left(f(\xi_n)\right)-	
    \ell_c(f) \phi\left(\frac{x}{\sqrt{n}}\right)\right| =0.
     \]
     Then the conclusions of Theorems \ref{thm1}, \ref{thm2} and \ref{thm4} remain true
     if  $h\in {\rm DRI}^+_c(\R)$ is replaced by $h\in B_c^+(\mathbb{Z})$, $\ell$ replaced by $\ell_c$ and $A$ replaced by $B\subset\mathbb{Z}$.
\end{remark}

\begin{lemma}\label{lemma1}
	Assume that {\bf(H2)}, {\bf(H3)} and {\bf(H4)} hold.
   For any $h\in {\rm DRI}^+_c(\R)$, it holds that
   	\begin{align}
   	    \lim_{t\to\infty} \sup_{x\in \R} \left|   \sqrt{t}\mathbf{E}_x\left(h(\xi_t)\right)-    \ell(h)
	    \phi\left(\frac{x}{\sqrt{t}}\right) \right| =0.
   \end{align}
\end{lemma}
\textbf{Proof: }
We write $t>1$ as $t=[t]+\gamma$ with $\gamma\in [0,1)$.
Combining the Markov property and the inequality $\sqrt{t}-\sqrt{[t]}= \frac{t-[t]}{\sqrt{t}+\sqrt{[t]}}\leq \frac{1}{\sqrt{t}}$,
we see that for any $x\in \R$,
\begin{align}
	&\left|   \sqrt{t}\mathbf{E}_x\left(h(\xi_t)\right)-	\ell(h)
	\phi\left(\frac{x}{\sqrt{t}}\right)\right| \nonumber\\
		&\leq (\sqrt{t}-\sqrt{[t]})\Vert h\Vert_\infty + \mathbf{E}_x\left( \left|   \sqrt{[t]}\mathbf{E}_{\xi_\gamma}(h(\xi_{[t]}))-\ell(h)\phi\left(\frac{\xi_\gamma}{\sqrt{[t]}}\right)  \right|\right)+\ell(h)\left|\mathbf{E}_x\left( \phi\left(\frac{\xi_\gamma}{\sqrt{[t]}}\right) \right)  - \phi\left(\frac{x}{\sqrt{t}}\right)\right|\nonumber\\
		& 	\leq \frac{\Vert h\Vert_\infty }{\sqrt{t}}
	+  \sup_{z\in \R} \left|   \sqrt{[t]}\mathbf{E}_z\left(h(\xi_{[t]})\right)-\ell(h)\phi\left(\frac{z}{\sqrt{[t]}}\right) \right| + \ell(h) \sup_{z\in\R} \left|\mathbf{E}_z\left( \phi\left(\frac{\xi_\gamma}{\sqrt{[t]}}\right) \right)  - \phi\left(\frac{z}{\sqrt{t}}\right)\right|.
\end{align}
The first term on the right-hand side tends to 0 as $t\to\infty$. By Lemma \ref{LLT-functional-xi}, the second term also tends to 0 as $t\to\infty$. Therefore, it remains to prove that
\begin{align}\label{Goal-1}
	\lim_{t\to\infty} \sup_{x\in\R}
		\left|\mathbf{E}_x\left( \phi\left(\frac{\xi_\gamma}{\sqrt{[t]}}\right) \right)  - \phi\left(\frac{x}{\sqrt{t}}\right)\right|
	= 	\lim_{t\to\infty} \sup_{x\in\R}
		\left|\mathbf{E}_0\left( \phi\left(\frac{\xi_\gamma +x}{\sqrt{[t]}}\right) \right)  - \phi\left(\frac{x}{\sqrt{t}}\right)\right|	=0.
\end{align}
Note that for $|x|> t^{2/3}$,
\begin{align}\label{Proof-Goal-1-1}
		& \left|\mathbf{E}_0\left( \phi\left(\frac{\xi_\gamma +x}{\sqrt{[t]}}\right) \right)  - \phi\left(\frac{x}{\sqrt{t}}\right)\right|\leq \mathbf{E}_0\left( \phi\left(\frac{\xi_\gamma +x}{\sqrt{[t]}}\right) \right)   +   \phi\left(\frac{x}{\sqrt{t}}\right)\nonumber\\
	&  \leq \mathbf{P}_0\left(\sup_{s\leq 1}|\xi_s| > \sqrt{t}\right) + \phi\left(\frac{t^{2/3}- \sqrt{t}}{\sqrt{[t]}}\right)+ \phi\left(t^{1/6}\right)\stackrel{t\to\infty}{\longrightarrow}0.
\end{align}
For $|x| \leq t^{2/3}$, by the inequality $|e^{-x^2}-e^{-y^2}| \leq |x^2 -y^2|$, we have
\begin{align}\label{Proof-Goal-1-2}
		& \left|\mathbf{E}_0\left( \phi\left(\frac{\xi_\gamma +x}{\sqrt{[t]}}\right) \right)  - \phi\left(\frac{x}{\sqrt{t}}\right)\right|\leq  \frac{1}{\sqrt{2\pi}} \mathbf{E}_0\left(\left| \frac{x^2-(\xi_\gamma+x)^2}{2[t]}\right|\right)+ \frac{x^2}{\sqrt{2\pi}}\left(\frac{1}{2[t]}-\frac{1}{2t}\right)\nonumber\\
	& \leq \frac{1}{2[t]\sqrt{2\pi}} \mathbf{E}_0\left(\xi_\gamma^2 + 2 |\xi_\gamma| |x|\right) + \frac{t^{1/3}}{2\sqrt{2\pi}[t]}\leq  \frac{1}{2[t]\sqrt{2\pi}} \left(\gamma + 2 \sqrt{\gamma} t^{2/3}\right)
	+\frac{t^{1/3}}{2\sqrt{2\pi}[t]}\nonumber\\
	&\leq   \frac{1}{2[t]\sqrt{2\pi}} \left(1 + 2  t^{2/3}\right)+\frac{t^{1/3}}{2\sqrt{2\pi}[t]} \stackrel{t\to\infty}{\longrightarrow}0.
\end{align}
Combining \eqref{Proof-Goal-1-1} and \eqref{Proof-Goal-1-2}, we get \eqref{Goal-1}. The proof is complete.
\hfill$\Box$
\bigskip

For $h\in {\rm DRI}^+_c(\R)$, define
\begin{align}\label{def-epsilon(t)}
	\epsilon_t(h):= \sup_{x\in\R}\left| \sqrt{t}\mathbf{E}_x\left(h(\xi_t)\right) -	\ell(h)
	\phi\left(\frac{x}{\sqrt{t}}\right) \right|\quad \mbox{and} \quad \widetilde{\epsilon}_t(h):= \sup_{q>t}\epsilon_q(h).
\end{align}
By the definition we easily see that, for any $t>0$, $\epsilon_t(h)\leq \sqrt{t} \Vert h\Vert_\infty +\frac{\ell(h)}{\sqrt{2\pi}}$.
 By Lemma \ref{lemma1}, we have $\sup_{t>0} \epsilon_t(h) + \sup_{t>0} \widetilde{\epsilon}_t(h)<\infty$ and
 \begin{equation}\label{limit0}
 \lim_{t\to\infty} \widetilde{\epsilon}_t(h)=\lim_{t\to\infty} \epsilon_t(h)=0.
 \end{equation}
Since $\ell(h)\phi(xt^{-1/2})\leq \frac{\ell(h)}{\sqrt{2\pi}}$ for any $x\in\R$, we have that, for any
$h\in {\rm DRI}_c^+(\R)$
and $g\in B_{Lip}^+(\R)$,
\begin{align}\label{uniformly-bounded-2}
	C(g, h)	:= \Vert g\Vert_\infty+ \sup_{x\in \R, t>0} \sqrt{t}\mathbf{E}_x(h(\xi_t))<\infty.
\end{align}

\begin{lemma}\label{lemma:small-effort}
	Assume that {\bf(H2)}, {\bf(H3)} and {\bf(H4)} hold.
		Let $h\in {\rm DRI}^+_c(\R)$
	and $g\in B_{Lip}^+(\R)$.
	\begin{enumerate}
		\item[(i)] For any $t, r>0$ and $y,z\in \R$, we have
		\begin{align}
				\sqrt{t}\left| \mathbf{E}_{\sqrt{t}y}(h(\xi_{tr}))-\mathbf{E}_{\sqrt{t}z}(h(\xi_{tr})) \right| \lesssim  \frac{\epsilon_{tr}(h)}{\sqrt{r}} +\frac{\sqrt{|y-z|}}{r^{3/4}}
		\end{align}
			and
		\begin{align}
			\left|\mathbf{E}_{\sqrt{t}y}\left(g\left(\frac{\xi_{tr}}{\sqrt{t}}\right)\right) -\mathbf{E}_{\sqrt{t}z}\left(g\left(\frac{\xi_{tr}}{\sqrt{t}}\right)\right)\right| \lesssim \left|y-z\right|.
		\end{align}
		\item[(ii)] For any $t>0, 0<r<s$ with $s-r\in (0,1) $ and $y\in \R$,  we have
		\begin{align}
				\sqrt{t}\left| \mathbf{E}_{\sqrt{t}y}(h(\xi_{tr}))-\mathbf{E}_{\sqrt{t}y}(h(\xi_{ts}) ) \right| \leq G_{h}(r,s; t),
		\end{align}
		where
		\begin{align}
			&G_{h}(r,s; t)\nonumber\\
			&:= \frac{\epsilon_{tr}(h)}{\sqrt{r}}+ \frac{\epsilon_{ts}(h)}{\sqrt{s}} + 	\ell(h)	\left(\frac{1}{\sqrt{r}}-\frac{1}{\sqrt{s}} + \frac{1}{\sqrt{s}}\left(\sqrt{s-r}+ 1-\exp\left\{-\frac{\sqrt{s-r}}{r}\right\} \right) \right).
		\end{align}
				Furthermore, for any $t>0, 0<r<s$  and $y\in \R$,  it holds that
		\begin{align}
				\left|\mathbf{E}_{\sqrt{t}y}\left(g\left(\frac{\xi_{tr}}{\sqrt{t}}\right)\right) -\mathbf{E}_{\sqrt{t}y}\left(g\left(\frac{\xi_{ts}}{\sqrt{t}}\right)\right)\right| \lesssim \sqrt{s-r}.
		\end{align}
	\end{enumerate}
\end{lemma}
\textbf{Proof: } (i)
The second inequality follows easily from
\begin{align}
	& 	\left|\mathbf{E}_{\sqrt{t}y}\left(g\left(\frac{\xi_{tr}}{\sqrt{t}}\right)\right) -\mathbf{E}_{\sqrt{t}z}\left(g\left(\frac{\xi_{tr}}{\sqrt{t}}\right)\right)\right|  \leq \mathbf{E}_0\left( \left| g\left(\frac{\xi_{tr}}{\sqrt{t}} +y\right) -g\left(\frac{\xi_{tr}}{\sqrt{t}} +z\right) \right|\right)\nonumber\\
	& \leq \mbox{Lip}(g)|y-z|\lesssim |y-z|.
\end{align}
Now we prove the first inequality.
By the definition of $\epsilon_t(h),$ we have
\begin{align}\label{step_33}
		\sqrt{t}\left| \mathbf{E}_{\sqrt{t}y}(h(\xi_{tr}))-\mathbf{E}_{\sqrt{t}z}(h(\xi_{tr})) \right| \leq \frac{2\epsilon_{tr}(h)}{\sqrt{r}}+ \frac{
			\ell(h)
		}{\sqrt{r}}\left|\phi\left(\frac{y}{\sqrt{r}}\right) -\phi\left(\frac{z}{\sqrt{r}}\right)  \right|.
\end{align}
Applying Lemma \ref{lemma2} (i) and \eqref{step_33}, we get the first inequality in  (i).

 (ii)  By H\"{o}lder's inequality,  we have
  \begin{align}
  	  &	\left|\mathbf{E}_{\sqrt{t}y}\left(g\left(\frac{\xi_{tr}}{\sqrt{t}}\right)\right) -\mathbf{E}_{\sqrt{t}y}\left(g\left(\frac{\xi_{ts}}{\sqrt{t}}\right)\right)\right| \leq \mathbf{E}_{\sqrt{t}y}\left(\left|  g\left(\frac{\xi_{tr}}{\sqrt{t}}\right) - g\left(\frac{\xi_{ts}}{\sqrt{t}}\right)  \right|\right)\nonumber\\
  	  & \leq \mbox{Lip}(g)  \mathbf{E}_{\sqrt{t}y}\left(\left|  \frac{\xi_{tr}}{\sqrt{t}} - \frac{\xi_{ts}}{\sqrt{t}}  \right|\right)\lesssim  \sqrt{ \mathbf{E}_{\sqrt{t}y}\left(\left|  \frac{\xi_{tr}}{\sqrt{t}} - \frac{\xi_{ts}}{\sqrt{t}}  \right|^2\right)} = \sqrt{s-r}.
  \end{align}
  The second inequality follows immediately.
  Now we prove the first inequality of (ii). Combining Lemma \ref{lemma2} (ii) and \eqref{step_33},
\begin{align}
	& 	\sqrt{t}\left| \mathbf{E}_{\sqrt{t}y}(h(\xi_{tr}))-\mathbf{E}_{\sqrt{t}y}(h(\xi_{ts}) ) \right|  \leq \frac{\epsilon_{tr}(h)}{\sqrt{r}}+ \frac{\epsilon_{ts}(h)}{\sqrt{s}} +
		\ell(h)
	 \left|\frac{1}{\sqrt{r}}\phi\left(\frac{y}{\sqrt{r}}\right)-\frac{1}{\sqrt{s}}\phi\left(\frac{y}{\sqrt{s}}\right)  \right|\nonumber\\
	& \leq \frac{\epsilon_{tr}(h)}{\sqrt{r}}+ \frac{\epsilon_{ts}(h)}{\sqrt{s}} +
	\ell(h)
	\left(\frac{1}{\sqrt{r}}-\frac{1}{\sqrt{s}} + \frac{1}{\sqrt{s}}\left(\sqrt{s-r}+ 1-\exp\left\{-\frac{\sqrt{s-r}}{r}\right\} \right)\right),
\end{align}
which implies the desired result.
\hfill$\Box$
\bigskip

For $x\in \R$, define
\begin{align}
	\tau_x^+:=\inf\{t>0 : \xi_t\geq x\},\quad \tau_x^-:=\inf\{t>0 : \xi_t\leq x\}.
\end{align}
The following result on the overshoot of $\xi$ is proved
in  \cite[Lemma 2.1]{HJRS}.

\begin{lemma}\label{Moment-of-overshoot}
	Assume that {\bf(H2)} holds.

	(i) If
	 $\mathbf{E}_0\left( ((-\xi_1)\vee 0)^\lambda \right)<\infty$ for some $\lambda>2$, then
	\begin{align}
		\sup_{x>0} \mathbf{E}_x\Big(\left|\xi_{\tau_0^-}\right|^{\lambda-2} \Big)<\infty.
	\end{align}
	
	(ii) If
	 $\mathbf{E}_0\left( (\xi_1\vee 0)^\lambda \right)<\infty$ for some $\lambda>2$, then
	\begin{align}
		\sup_{x>0} \mathbf{E}_{-x}\Big(\xi_{\tau_0^+}^{\lambda-2} \Big)<\infty.
	\end{align}
\end{lemma}

\subsection{Evolution equation for $(Z_t)_{t\ge0}$}
In this section, we always assume that {\bf(H1)}--{\bf(H4)} hold.

For any $f\in B^+(\R)$, define
\begin{align}
		    v_f(t,y) &:=1- \E_y\left(\exp\left\{-\int_\R f(x)Z_t(\mathrm{d}x) \right\}\right).
\end{align}
The next lemma gives an integral equation for $v_f$.
Using \cite[Lemma 4.1]{E.B1.}, its proof is standard and  similar to that in \cite[Lemma 2.1]{HRS}, and so we omit it.
Define
\begin{align}
	\psi(v):=\beta\left(\sum_{k=0}^\infty p_k(1-v)^k-(1-v)\right),\quad v\in [0,1].
\end{align}
Since $\sum_{k=0}^\infty kp_k=1$, by Jensen's inequality, we have that
\begin{align}\label{Property-phi-non-negative}
	\psi(v)\geq \beta\left( (1-v)^{\sum_{k=0}^\infty kp_k} -(1-v)\right)=0.
\end{align}

\begin{lemma}\label{Integration-equation}
	For any $t>0, y\in\R$, $v_f(t,y)$ solves 	the equation
	\begin{align}
		v_f(t,y)=\mathbf{E}_y\left(1-e^{-f(\xi_t)}\right) -\mathbf{E}_y\left(\int_0^t \psi\left(v_f(t-s,\xi_s)\right)\mathrm{d}s\right).
	\end{align}
\end{lemma}

For any $t ,r>0$, we define
\begin{align}\label{Def-of-parameters-2}
\psi^{(t)}(v) := t^{\frac{\alpha}{\alpha-1}}\psi\left(v t^{-\frac{1}{\alpha-1}}\right),\qquad \xi_r^{(t)}:=\frac{\xi_{tr}}{\sqrt{t}}.
\end{align}
For any $h\in \mathcal{C}_c^+(\R), g\in B_{Lip}^+(\R)$, $t ,r>0$ and $y\in \R$, we define
\begin{align}\label{Def-of-parameters}
	f^{(t)}(\cdot)& := \frac{1}{t^{\frac{1}{\alpha-1}-\frac{1}{2}}}h(\cdot) + \frac{1}{t^{\frac{1}{\alpha-1}}}g\left(\frac{\cdot}{\sqrt{t}}\right),\quad
		v_{g, h}^{(t)}(r,y):= t^{\frac{1}{\alpha-1}}
	v_{f^{(t)}}(tr, \sqrt{t}y).
\end{align}
With a slight abuse of notation, we will
use the same notation $\mathbf{P}_y$ to denote the law of $\xi_r^{(t)}$ starting from
 $\xi_0^{(t)}=y$, i.e.,
\begin{align}\label{scaling}
(\xi_r^{(t)},\mathbf{P}_y)
:=
\left(\frac{1}{\sqrt{t}}\xi_{tr}, \mathbf{P}_{\sqrt{t}y}\right).
\end{align}
It is well-known (for example, see
\cite[Lemma 2.14(ii)]{HRS})
that, under {\bf(H1)}, for any $K>0$,
\begin{align}\label{Properties-of-phi}
\lim_{t\to\infty} \psi^{(t)}(v)= \mathcal{C}(\alpha) v^\alpha\ \mbox{ uniformly for $v\in [0,K]$.}
\end{align}

\begin{lemma}\label{lemma3}
	There exists a constant $C_\psi \in (0,\infty)$ such that
	\begin{align}
		\left| \psi^{(t)}(u)- \psi^{(t)}(v) \right|\leq C_\psi (u^{\alpha-1}+v^{\alpha-1})|u-v|,\quad \forall u,v\in [0,t^{\frac{1}{\alpha-1}}]
		,\quad \forall t>0.
	\end{align}
	In particular, we have
	\begin{align}
		 \psi^{(t)}(v) \leq C_\psi v^{\alpha},\quad \forall v\in [0,t^{\frac{1}{\alpha-1}}]
		 ,\quad\forall t>0.
	\end{align}
\end{lemma}
\textbf{Proof: } We first prove that there exists some constant $C_\psi$ such that
\begin{align}\label{step_9}
	\left|\psi'(v)\right|\leq C_\psi v^{\alpha-1},\quad \forall v\in [0,1].
\end{align}
First, using $\sum^\infty_{k=1}kp_k=1$, we have
\begin{align}\label{step_65}
	& \left|\psi'(v)\right|= \beta \left(1-\sum_{k=1}^\infty kp_k (1-v)^{k-1}\right)= \beta v\sum_{k=2}^\infty kp_k \left(\sum_{j=0}^{k-2}(1-v)^j\right)\nonumber\\
	& = \beta v\sum_{j=0}^\infty (1-v)^j \sum_{k=j+2}^\infty kp_k.
\end{align}
Under {\bf(H1)}, we have $\sum_{k=n}^\infty p_k \lesssim n^{-\alpha}$ for all $n\geq 2$. Thus, for all $j\geq 0$,
\begin{align}\label{step_66}
	& \sum_{k=j+2}^\infty kp_k = (j+1) \sum_{k=j+2}^\infty p_k + \sum_{k=j+2}^\infty \sum_{n=k}^\infty p_n \nonumber\\& \lesssim \frac{j+1}{(j+2)^\alpha}+ \sum_{k=j+2}^\infty \frac{1}{k^\alpha}\lesssim  \sum_{k=j+2}^\infty \frac{1}{k^\alpha}.
\end{align}
Combining \eqref{step_65} and   \eqref{step_66}, we conclude that for all $v\in [0,1]$,
\begin{align}
	& \left|\psi'(v)\right| \lesssim  v\sum_{j=0}^\infty (1-v)^j\sum_{k=j+2}^\infty \frac{1}{k^\alpha} =  \sum_{k=2}^\infty \frac{1}{k^\alpha} \left(1-(1-v)^{k-1}\right).
\end{align}
Together with the inequality $1-(1-v)^{k-1}\leq 1\land ((k-1)v)$, we obtain that for all $v\in [0,1]$,
\begin{align}
		& \left|\psi'(v)\right| \lesssim \sum_{k=2}^\infty \frac{1}{k^\alpha}
				(1\land ((k-1)v))
		\leq  \int_1^\infty \frac{1}{x^\alpha} (1\land (xv))\mathrm{d}x =   v\int_1^{1/v} \frac{1}{x^{\alpha-1}}\mathrm{d}x + \int_{1/v}^\infty \frac{1}{x^\alpha}\mathrm{d}x,
\end{align}
which implies \eqref{step_9}.

Now we assume that $u<v$. Then there exists
 $\eta \in [u,v]$ such that
\begin{align}
	& \left| \psi^{(t)}(u)- \psi^{(t)}(v) \right| = t^{\frac{\alpha}{\alpha-1}}\left|\psi(ut^{-\frac{1}{\alpha-1}})-\psi(vt^{-\frac{1}{\alpha-1}})\right| \nonumber\\
	& = t\left|u-v\right|
		    \left|\psi'(\eta t^{-\frac{1}{\alpha-1}})\right|\leq C_\psi \eta^{\alpha-1}|u-v|
	\leq C_\psi v^{\alpha-1}|u-v|,
\end{align}
where in the second to last inequality, we used \eqref{step_9}. The proof is complete.
\hfill$\Box$

\begin{lemma}\label{Scaled-evolution-equation}
	For any
	$h\in {\rm DRI}^+_c(\R)$
	and $g\in B_{Lip}^+(\R)$, it holds that
	\begin{align}
		& v_{g, h}^{(t)}(r,y)	= t^{\frac{1}{\alpha-1}}\mathbf{E}_{y}\left(1-\exp\left\{-\frac{1}{t^{\frac{1}{\alpha-1}-\frac{1}{2}}}h(\sqrt{t}\xi_{r}^{(t)}) -\frac{1}{t^{\frac{1}{\alpha-1}}}g(\xi_r^{(t)})\right\}\right) \nonumber\\
		&\quad\quad - \mathbf{E}_{y}\left(\int_0^r \psi^{(t)} \left(v_{g, h}^{(t)}	(r-s,\xi_s^{(t)})\right)\mathrm{d}s\right).
	\end{align}
\end{lemma}
\textbf{Proof:}
Combining \eqref{Def-of-parameters-2}, \eqref{Def-of-parameters} and Lemma \ref{Integration-equation}, we get that
\begin{align}
		& v_{g, h}^{(t)}(r,y)
	=t^{\frac{1}{\alpha-1}}\mathbf{E}_{y}\left(1-\exp\left\{-\frac{1}{t^{\frac{1}{\alpha-1}-\frac{1}{2}}}h(\sqrt{t}\xi_{r}^{(t)}) -\frac{1}{t^{\frac{1}{\alpha-1}}}g(\xi_r^{(t)})\right\}\right) \nonumber\\
	&\quad -  t^{\frac{1}{\alpha-1}}\mathbf{E}_{\sqrt{t}y}\left(\int_0^{tr} \psi\left(
		v_{g, h}^{(t)}
	(tr-s,\xi_s)\right)\mathrm{d}s\right)\nonumber\\
	& = t^{\frac{1}{\alpha-1}}\mathbf{E}_{y}\left(1-\exp\left\{-f^{(t)}(\sqrt{t}\xi_{r}^{(t)})\right\}\right)  -  t^{\frac{1}{\alpha-1}}\mathbf{E}_{y}\left(\int_0^{tr} \psi\left(
		v_{g, h}^{(t)}
	(tr-s, \sqrt{t}\xi_{s/t}^{(t)})\right)\mathrm{d}s\right)\nonumber\\
	& = t^{\frac{1}{\alpha-1}}\mathbf{E}_{y}\left(1-\exp\left\{-f^{(t)}(\sqrt{t}\xi_{r}^{(t)})\right\}\right) -  t^{\frac{\alpha}{\alpha-1}}\mathbf{E}_{y}\left(\int_0^{r} \psi\left(
		v_{g, h}^{(t)}
	(tr-ts, \sqrt{t}\xi_{s}^{(t)})\right)\mathrm{d}s\right)\nonumber\\
	& = t^{\frac{1}{\alpha-1}}\mathbf{E}_{y}\left(1-\exp\left\{-f^{(t)}(\sqrt{t}\xi_{r}^{(t)})\right\}\right) -  t^{\frac{\alpha}{\alpha-1}}\mathbf{E}_{y}\left(\int_0^{r} \psi\left( t^{-\frac{1}{\alpha-1}}
		v_{g, h}^{(t)}
	(r-s, \xi_{s}^{(t)})\right)\mathrm{d}s\right).
\end{align}
The desired result now follows immediately from the definition of $\psi^{(t)}$.
\hfill$\Box$

\subsection{Initial trace theory}

For any open set $U\subset \R$, we denote by $C_c^+(U)$ the family of non-negative continuous functions with compact support in $U$.
Denote by $\mathcal{B}_{reg}^+(\R)$  the space of positive outer regular Borel measures.
Suppose that $\Lambda\subset \R$ is a closed set and that $\eta$ is a non-negative Radon measure on $\Lambda^c$.  By
 \cite[pp.1452--1453]{MV1999},
 the pair $(\Lambda, \eta)$ can be represented by the following measure $\gamma_{(\Lambda,\eta)}\in \mathcal{B}_{reg}^+(\R)$:
\begin{equation}\label{Def-Regular-measure-gamma}
	\gamma_{(\Lambda,\eta)}(B):=
\left\{\begin{array}{ll}
       \infty,\quad
			&B\cap \Lambda \neq \emptyset,\\
			\eta(B),\quad &
			B\cap \Lambda  = \emptyset.
		\end{array}\right.
\end{equation}
Define  the set of regular points of $\gamma_{(\Lambda,\eta)}$ by
\[
\mathcal{R}_{(\Lambda, \eta)}:= \left\{x\in \R: \gamma_{(\Lambda,\eta)} ((x-z,x+z))=\infty,\quad \forall z>0\right\}^c.
\]

For any closed set $\hat{\Lambda}\subset \R$ and non-negative Radon measure $\hat{\eta}$ on $\hat{\Lambda}^c$, consider the problem
 \begin{equation}\label{PDE-2}
	\left\{\begin{array}{rl}
		&\frac{\partial}{\partial r} \hat{v}_{(\hat{\Lambda},\hat{\eta})}^X(r,y) = \frac{\partial^2}{\partial y^2} \hat{v}_{(\hat{\Lambda},\hat{\eta})}^X(r,y)-  \left(\hat{v}_{(\hat{\Lambda},\hat{\eta})}^X(r,y)\right)^\alpha, \, (r,y)\in (0,\infty)\times \R, \\
		& \left\{x \in \R: \forall z>0, \ \lim_{r\downarrow0} \int_{x-z}^{x+z}\hat{v}_{(\hat{\Lambda},\hat{\eta})}^X(r,y) \mathrm{d} y=\infty \right\}=\hat{\Lambda},\\
		& \forall f\in C_c^+(\hat{\Lambda}^c),\  \lim_{r\downarrow 0} \int f(y)   \hat{v}_{(\hat{\Lambda},\hat{\eta})}^X(r,y) \mathrm{d} y = \int f(y) \hat{\eta}(\mathrm{d}y).
	\end{array}\right.
\end{equation}
Define $\Lambda:= \{x/\sqrt{2}: x\in \hat{\Lambda}\}$ and let $\eta$ be the Radon measure on $\Lambda^c$ such that
\[
\int_{\Lambda^c} f(y)\eta(\mathrm{d}y): =\frac{1}{\sqrt{2}\mathcal{C}(\alpha)^{\frac{1}{\alpha-1}}}\int_{\hat{\Lambda}^c} f(\sqrt{2}y)\hat{\eta}(\mathrm{d} y).
\]
Consider the problem
\begin{equation}\label{PDE}
	\left\{\begin{array}{rl}
		&\frac{\partial}{\partial r} v_{(\Lambda,\eta)}^X(r,y) = \frac{1}{2} \frac{\partial^2}{\partial y^2} v_{(\Lambda,\eta)}^X(r,y)- \varphi\left(v_{(\Lambda,\eta)}^X(r,y)\right), \, (r,y)\in (0,\infty)\times \R, \\
		& \left\{x \in \R: \forall z>0, \ \lim_{r\downarrow0} \int_{x-z}^{x+z}v_{(\Lambda,\eta)}^X(r,y) \mathrm{d} y=\infty \right\}=\Lambda,\\
		& \forall f\in C_c^+(\Lambda^c),\  \lim_{r\downarrow 0} \int f(y)   v_{(\Lambda,\eta)}^X(r,y) \mathrm{d} y = \int f(y)\eta(\mathrm{d}y).
	\end{array}\right.
\end{equation}
It is easy to check that
\[
v_{(\Lambda,\eta)}^X(r,y) = \frac{1}{\mathcal{C}(\alpha)^{\frac{1}{\alpha-1}}} \hat{v}_{(\hat{\Lambda},\hat{\eta})}^X\left(r,\frac{y}{\sqrt{2}}
\right)
\]
is a one-to-one correspondence between the positive solutions of \eqref{PDE-2} and the positive solutions of \eqref{PDE}. According to \cite[Theorem 3.5]{MV1999},
\eqref{PDE-2} has a unique positive solution
$\hat{v}_{(\hat{\Lambda},\hat{\eta})}^X(r,y)$.
Consequently,
the function $v^X_{(\Lambda, \eta)}(r,y)$ defined above is the unique solution of \eqref{PDE}. We call
$(\Lambda, \eta)$  the initial trace of the solution $v_{(\Lambda,\eta)}^X$.

 In this section, we give a  probabilistic representation of the solution $v^X_{(\Lambda, \eta)}$.
 To avoid too much measure theoretic details, we only deal with the case when $\Lambda$ is a bounded closed interval.
 In the special case $\varphi(\lambda)=\frac{1}{2}\lambda^2$, a probabilistic representation via Brownian snake was given
 by  Le Gall \cite[Theorem 4]{LeGall1996}.

 Recall that $X$ is a  critical super-Brownian motion with branching mechanism $\varphi$ given in \eqref{Stable-Branching-mechanism} and $\{Y_t(x): t>0, x\in \R\}$ is the density process of $X$ which, for all $y\in \R$, is $\P_{\delta_y}$-almost surely continuous
  in $x$ for all $t>0$.
Let  $v_{(\Lambda,\eta)}^X$ is the solution of the PDE problem \eqref{PDE}.

\begin{prop}\label{Probabilistic-representation}
	Suppose that $\Lambda=[a,b]\subset \R$ is a bounded closed interval and $\eta$ is a Radon measure on $\Lambda^c$.
  Then for any $r>0, y\in\R$,
	\begin{align}\label{Identity}
		 v_{(\Lambda,\eta)}^X(r,y)= -\log \E_{\delta_y}\left(1_{\{Y_r(x)=0,\ \forall x\in \Lambda\}} e^{-\int Y_r(z) \eta(\mathrm{d}z)}\right).
	\end{align}
\end{prop}
Before presenting the proof, we first recall the notion of $m$-weak convergence from \cite[Definition 3.9]{MV1999}.

\begin{defn}\label{m-conv}
Let $(\Lambda_n, \eta_n)$ be a sequence of initial traces and $(\Lambda,\eta)$ be another initial trace.
We say that the measures $\gamma_{(\Lambda_n,\eta_n)}$ converge $m$-weakly to the measure $\gamma_{(\Lambda,\eta)}$ if the following
two conditions hold:
\begin{description}
\item[(i)] If $U\subset\R$ is an open set with $\gamma_{(\Lambda,\eta)}(U)=\infty$, then $\lim_{n\to\infty} \gamma_{(\Lambda_n, \eta_n)}(U)=\infty$.

\item[(ii)] For any compact set $K\subset \mathcal{R}_{(\Lambda, \eta)}$, the sequence of $\gamma_{(\Lambda_n,\eta_n)}(K)$ is eventually bounded, i.e., there exists $N\in \N$ and $C\in (0,\infty)$ such that $\gamma_{(\Lambda_n,\eta_n)}(K)\leq C$ for all $n\geq N$,
 and  for any $\phi\in C_c^+(\mathcal{R}_{(\Lambda, \eta)})$,
$\lim_{n\to\infty} \int \phi(x)\gamma_{(\Lambda_n,\eta_n)}(\mathrm{d}x)=\int \phi(x)\gamma_{(\Lambda,\eta)}(\mathrm{d}x)$.
\end{description}
\end{defn}

According to  \cite[Section 2.1]{BMS2024}, if
$\gamma_{(\Lambda_n, \eta_n)}$ converges $m$-weakly to $\gamma_{(\Lambda,\eta)}$,
then for any $r>0$ and $y\in\R$, $v_{(\Lambda_n, \eta_n)}^X(r,y)$ converges to $v_{(\Lambda,\eta)}^X(r,y)$ as $n\to \infty$.
Now we are ready to prove Proposition \ref{Probabilistic-representation}.

\bigskip

\textbf{Proof of Proposition \ref{Probabilistic-representation}: }
Step 1:
In this part we consider the case that $\Lambda=\emptyset$. It is well-known that for any $r>0$, $Y_r(\cdot)$ is compactly supported (to see this, one can fix $x$ and $t$ and
let the constant $\Lambda$ in \cite[Lemma 4.3]{RSZ2021} tend to $\infty$).
Therefore, by the Markov property and the dominated convergence theorem,
\begin{align}
	& v_{(\emptyset,\eta)}^X(r,y)= \lim_{s\downarrow 0}	v_{(\emptyset,\eta)}^X(r+s,y)= -\lim_{s\downarrow 0}\log \E_{\delta_y} \left(\exp\left\{-\int v_{(\emptyset,\eta)}^X(s,z) Y_r(z)\mathrm{d}z \right\}\right)\nonumber\\ & = -\log \E_{\delta_y} \left(\lim_{s\downarrow 0} \exp\left\{-\int v_{(\emptyset,\eta)}^X(s,z) Y_r(z)\mathrm{d}z \right\}\right) = -\log \E_{\delta_y} \left( \exp\left\{-\int  Y_r(z)\eta(\mathrm{d}z) \right\}\right),
\end{align}
which implies \eqref{Identity} in the case $\Lambda=\emptyset$.

Step 2:
In this step we consider the case that $\Lambda\subset \R$ is a closed subset and that $\eta$ is a Radon measure on  $\Lambda^c$. Define $\eta_\Lambda(\mathrm{d}x)= 1_{\Lambda}\mathrm{d}x$ if
$\ell(\Lambda)=b-a\neq 0$
 and $\eta_\Lambda(\mathrm{d}x)= \delta_a(\mathrm{d}x)$ if $\Lambda=\{a\}$.
For each $n$, define
\begin{align}
	\Lambda_n:= \emptyset,\quad B_n:= \left\{y\in \R: \mbox{dist}(y, \Lambda) \leq \frac{1}{n}\right\}\quad\mbox{and}\quad \eta_n:= n\eta_\Lambda+ \eta|_{B_n^c},
\end{align}
then for any $n$, $\eta_n$ is a Radon measure on $\R$. By the result obtained in Step 1,
we have
\begin{align}\label{step_11}
	&v_{(\Lambda_n, \eta_n)}^X(r,y)= -\log \E_{\delta_y} \left( \exp\left\{-\int  Y_r(z)\eta_n(\mathrm{d}z) \right\}\right)\nonumber\\
	& = -\log \E_{\delta_y} \left( \exp\left\{-n \int_\Lambda  Y_r(z)\eta_\Lambda(\mathrm{d}z) - \int_{B_n^c} Y_r(z)\eta(\mathrm{d}z) \right\}\right).
\end{align}
Since $B_n^c \uparrow \Lambda^c$ as $n\to\infty$, combining the dominated convergence theorem and \eqref{step_11}, we see that
\begin{align}
	&\lim_{n\to\infty} v_{(\Lambda_n, \eta_n)}^X(r,y) = -\log \E_{\delta_y} \left(1_{\{\int_\Lambda Y_r(z)\eta_\Lambda(z)=0\}} \exp\left\{- \int_{\Lambda^c} Y_r(z)\eta(\mathrm{d}z) \right\}\right)\nonumber\\
	&= -\log \E_{\delta_y} \left(1_{\{Y_r(z)=0,\ \forall z\in \Lambda\}} \exp\left\{- \int Y_r(z)\eta(\mathrm{d}z) \right\}\right),
\end{align}
where in the last equality
we used the fact that the support of $\eta_\Lambda$ is equal to $\Lambda$ and that the support of $\eta$ is a subset of $\Lambda^c$.
Therefore,  to complete the proof, it suffices to show that $(\Lambda_n, \eta_n)$ converges $m$-weakly to $(\Lambda,\eta)$.
Now we check the conditions in Definition \ref{m-conv}.
For (i), suppose that $\gamma_{(\Lambda,\eta)}(U)=\infty$ for an open set $U\subset \R$. If $U\cap \Lambda\neq \emptyset$, then by the definition of $\eta_\Lambda$,  we can find a Borel set $B\subset U\cap \Lambda$
such that $\eta_\Lambda(B)>0$. In this case,
\[
\eta_n(U)\geq \eta_n(B)= n \eta_\Lambda(B)\stackrel{n\to\infty}{\longrightarrow}\infty.
\]
On the other hand, if $U\cap \Lambda=\emptyset$, then
\[
\eta(U)=\infty= \eta(U\cap \Lambda^c)= \lim_{n\to\infty} \eta(U\cap B_n^c) \leq \lim_{n\to\infty}  \eta_n(U),
\]
as desired.  Now we check (ii).
Note that $\mathcal{R}_{(\Lambda,\eta)}= \Lambda^c$. For
any compact set $K\subset \Lambda^c$, $\eta_n(K)= \eta(K\cap B_n^c)\leq \eta(K)<\infty$ is bounded. Besides, for any $\phi\in C_c^+(\Lambda^c)$, by
the monotone convergence theorem,
\begin{align}
	\lim_{n\to\infty} \int \phi(x)\gamma_{(\Lambda_n,\eta_n)}(\mathrm{d}x)=\lim_{n\to\infty} \int_{B_n^c} \phi(x)\eta(\mathrm{d}x)=\int_{\Lambda^c} \phi(x)\eta(\mathrm{d}x)  =\int \phi(x)\gamma_{(\Lambda,\eta)}(\mathrm{d}x),
\end{align}
which implies (ii). This completes the proof of the proposition.

\hfill$\Box$

\begin{remark}\label{remark}
		We will need the following result later: for  any $r>0, y\in \R$,
	\begin{align}\label{step_26}
		& \lim_{\varepsilon\to 0} v_{([-\varepsilon, \varepsilon], 0)}^X(r,y) = -\lim_{\varepsilon\to 0} \log \E_{\delta_y}\left(Y_r(x)=0,\ \forall x\in [-\varepsilon, \varepsilon]\right)\nonumber\\
		& = v_{(\{0\}, 0)}^X(r,y) = -\log \E_{\delta_y}\left(Y_r(0)=0\right).
	\end{align}
	To prove \eqref{step_26}, we only need to show that $\gamma_{([-\varepsilon, \varepsilon], 0)}$ converges $m$-weakly to $\gamma_{(\{0\},0)}.$ Condition (i) of
Definition \ref{m-conv}
is easy to check since $\gamma_{([-\varepsilon, \varepsilon], 0)}(U)\geq \gamma_{(\{0\},0)}(U).$ For (ii), for conpact set $K\subset \mathcal{R}_{(\{0\}, 0)}= \R\setminus \{0\}$, let $\varepsilon_0>0$ sufficient small so that $K\subset [-\varepsilon_0, \varepsilon_0]^c$. Then $\gamma_{([-\varepsilon, \varepsilon], 0)}(K)=0$ when $\varepsilon<\varepsilon_0$. Furthermore, for any $\phi\in C_c^+(\R \setminus \{0\})$, suppose that the support of $\phi$ is a subset of $[-\varepsilon_0,\varepsilon_0]^c$, then for any $\varepsilon<\varepsilon_0$, it holds that $ \int \phi(x)\gamma_{([-\varepsilon,\varepsilon],0)}(\mathrm{d}x)=0=\int \phi(x)\gamma_{(\{0\},0)}(\mathrm{d}x)$. Hence \eqref{step_26} is true.
\end{remark}

\section{Proof of the main results}\label{S:3}

In this section, we always assume that {\bf(H1)}--{\bf(H4)} hold.

\begin{lemma}\label{lemma4}
		Let $h\in {\rm DRI}^+_c(\R)$
		and $g\in B_{Lip}^+(\R)$, and let	$v_{g, h}^{(t)}$
		be given in \eqref{Def-of-parameters}.
       Suppose $r>0$.
      Then there exists a constant
       $C_1=C_1(g, h)$
      such that for any $t>1, y\in \R$ and $s\in [0,r]$,
		\begin{align}\label{step_62}
			\mathbf{E}_y \left(\left(	v_{g, h}^{(t)}	(r-s,\xi_s^{(t)})\right)^{\alpha-1}\right)   \leq \left(	\mathbf{E}_y \left(	v_{g, h}^{(t)}	(r-s,\xi_s^{(t)})\right)\right)^{\alpha-1}\leq
		 C_1
		 \left(\frac{1}{r^{(\alpha-1)/2}} \land t^{(\alpha-1)/2}\right)
		\end{align}
		and that
		\begin{align}\label{step_67}
			\mathbf{E}_y \left(\psi^{(t)}\left(	v_{g, h}^{(t)} (r-s,\xi_s^{(t)})\right)  \right)     \leq
				 \frac{C_1}{(r-s)^{(\alpha-1)/2}\sqrt{r}}.
		\end{align}
\end{lemma}
\textbf{Proof: }
The first inequality of \eqref{step_62} follows directly from Jensen's inequality. Now we prove the second inequality of \eqref{step_62}.
Combining Lemma \ref{Scaled-evolution-equation},  \eqref{uniformly-bounded-2}, \eqref{scaling}
and the fact that $1-e^{-|x|}\leq |x|$,
we get that for all $t>1, y\in \R$ and $r>0$,
\begin{align}\label{upper-bound-of-v-theta-A}
		&v_{g, h}^{(t)}
	(r,y)\leq  \sqrt{t} \mathbf{E}_y\left( h(\sqrt{t}\xi_r^{(t)})\right)  + \mathbf{E}_y\left(g(\xi_r^{(t)})\right)\nonumber\\
		&\leq \left(
		\frac{C(g,h)}
		{\sqrt{r}}\right)\land (\sqrt{t}\Vert h\Vert_\infty+\Vert g\Vert_\infty)
	\lesssim \frac{1}{\sqrt{r}} \land \sqrt{t}.
\end{align}
Therefore, combining \eqref{upper-bound-of-v-theta-A} and the Markov  property, for any $t>1, y\in \R, s\in [0,r]$,
\begin{align}
	& \mathbf{E}_y \left(	v_{g, h}^{(t)}	(r-s,\xi_s^{(t)})\right)\leq \mathbf{E}_y \left(\sqrt{t} \mathbf{E}_{\xi_s^{(t)}}\left( h(\sqrt{t}\xi_{r-s}^{(t)})\right)  + \mathbf{E}_{\xi_s^{(t)}}\left(g(\xi_{r-s}^{(t)})\right) \right) \nonumber\\
	& =   \sqrt{t} \mathbf{E}_y\left( h(\sqrt{t}\xi_r^{(t)})\right)  + \mathbf{E}_y\left(g(\xi_r^{(t)})\right)\lesssim \frac{1}{\sqrt{r}} \land \sqrt{t},
\end{align}
 which implies \eqref{step_62}.
For \eqref{step_67}, combining Lemma \ref{lemma3}, \eqref{step_62} and \eqref{upper-bound-of-v-theta-A},
\begin{align}\label{step_52}
	& \mathbf{E}_y \left(\psi^{(t)}\left(v_{g, h}^{(t)}
	(r-s,\xi_s^{(t)})\right)\right) \lesssim  \mathbf{E}_y \left(\left(	v_{g, h}^{(t)}
	(r-s,\xi_s^{(t)})\right) ^\alpha \right)\nonumber\\
	& \lesssim \frac{1}{(r-s)^{(\alpha-1)/2}} \mathbf{E}_y \left(	v_{g, h}^{(t)}
	(r-s,\xi_s^{(t)})\right)\lesssim \frac{1}{(r-s)^{(\alpha-1)/2} \sqrt{r}}.
\end{align}
The proof is now complete.
\hfill$\Box$

Recall the definition of $\widetilde{\epsilon}_t(h)$ defined in  \eqref{def-epsilon(t)}.

\begin{prop}\label{prop1}
	Assume
	$h\in {\rm DRI}^+_c(\R)$,
	$g\in B_{Lip}^+(\R)$ and $T>0$. Let $v_{g, h}^{(t)}$ be defined as in \eqref{Def-of-parameters}.
	
    (i) There exists a constant $N_1=N_1(g, h,  T, \psi)>0$  such that for all $t>1, r\in (0, T]$ and $y,z\in \R$ with $|y-z|<1$,
    \begin{align}\label{step_34}
    \left|   v_{g, h}^{(t)} (r,y)- v_{g, h}^{(t)} (r,z) \right|\leq \frac{N_1}{r^{3/4}}
    \left( \widetilde{\epsilon}_{\sqrt{t}r}(h)+ \sqrt{|y-z|}+ \frac{1}{t^{1/4}} \right).
    \end{align}
    (ii) There exists a constant  $N_2=N_2(g, h, T, \psi)>0$ such that for all $t>1$, $r\in (0,T]$, $q\in (0,1)$ and $y\in\R$,
    \begin{align}\label{step_35}
   	\left|	v_{g, h}^{(t)}	(r,y)-v_{g, h}^{(t)} (r+q,y) \right|\leq	\frac{N_2}{r^{3/2}}
	\left(\widetilde{\epsilon}_{\sqrt{t}r}(h)+  q^{1/8}+\frac{1}{t^{1/4}}\right).
    \end{align}
\end{prop}
\textbf{Proof:} (i)   Without loss of generality, we assume that $y<z$.
Combining
\eqref{Lemma:inequality} with Lemmas \ref{lemma3} and  \ref{Scaled-evolution-equation},
we get that for all $t>1, y,z\in \R$ and $r\in (0, T]$,
\begin{align}\label{step_53}
& \left|v_{g, h}^{(t)}(r,y)-v_{g, h}^{(t)}(r,z) \right|\nonumber\\
\le &\frac{t^{\frac{1}{\alpha-1}}}{t^{\frac{2}{\alpha-1}-1}} \left(\mathbf{E}_y\left( h^2 (\sqrt{t}\xi_r^{(t)})\right) +\mathbf{E}_z\left( h^2 (\sqrt{t}\xi_r^{(t)})\right)  \right)+ \frac{1}{t^{\frac{1}{\alpha-1}}} \left(  \mathbf{E}_y\left(g^2(\xi_r^{(t)})  \right)+\mathbf{E}_z\left(g^2(\xi_r^{(t)})  \right) \right)\nonumber\\
&+\sqrt{t}\sup_{x\in\R} \left|\mathbf{E}_x(h(\sqrt{t}\xi_r^{(t)}))- \mathbf{E}_{x+y-z}(h(\sqrt{t}\xi_r^{(t)}))\right|
+ \sup_{x\in \R}  \left|\mathbf{E}_x(g(\xi_r^{(t)}))- \mathbf{E}_{x+y-z}(g(\xi_r^{(t)}))\right|\nonumber\\
&+C_\psi \int_0^{r} \left(\mathbf{E}_y\left( \left(v_{g, h}^{(t)}(r-s,\xi_s^{(t)})\right)^{\alpha-1} \right) +\mathbf{E}_z\left( \left(v_{g, h}^{(t)}(r-s,\xi_s^{(t)})\right)^{\alpha-1} \right)\right) \nonumber\\
&\quad\quad\times\sup_{x\in \R}\left|v_{g, h}^{(t)}(r-s,x)-v_{g, h}^{(t)}(r-s,x+z-y)\right| \mathrm{d}s.
\end{align}
Define
\[
F_{z-y}^f(r):=\sup_{x\in\R} \left| v_{g, h}^{(t)} (r, x)- v_{g, h}^{(t)}(r, x+z-y)\right|.
\]
Combining Lemma \ref{lemma:small-effort}(i),   \eqref{uniformly-bounded-2}, \eqref{scaling},
Lemma \ref{lemma4} and the fact that $\frac{1}{\alpha-1}\geq 1$,  we conclude from the above inequality that  for any $y,z\in \R, r\in (0, T]$ and $t>1$,
\begin{align}\label{step_54}
	& F_{z-y}^f(r) \lesssim  \frac{1}{\sqrt{tr}} \land 1 + \frac{1}{t}+\left(\frac{\epsilon_{tr}(h)}{\sqrt{r}}+ \frac{\sqrt{|y-z|}}{r^{3/4}} \right)\land \sqrt{t} +  |y-z|+   \frac{1}{r^{(\alpha-1)/2}}\int_{0}^{r} F_{z-y}^f(s)\mathrm{d}s .
\end{align}
Define $G_{z-y}^f(r):= r^{(\alpha-1)/2}F_{z-y}^f(r)$. Then there exists a constant $K_1=K_1(g, h, \psi, T) \in (0,\infty)$ such that for all $y,z\in \R$ with $|y-z|<1$, $t>1$ and $r\in (0,T]$,
\begin{align}
	G_{z-y}^f(r)&\leq K_1
	r^{(\alpha-1)/2} \left( \frac{1}{\sqrt{tr}} \land 1 + \frac{1}{t}+\left(\frac{\epsilon_{tr}(h)}{\sqrt{r}}+ \frac{\sqrt{|y-z|}}{r^{3/4}} \right)\land \sqrt{t} +  |y-z|\right) \nonumber\\
	&\quad +K_1	\int_{0}^{r} \frac{1}{s^{(\alpha-1)/2}}G_{z-y}^f(s)\mathrm{d}s \nonumber\\
	&=: \alpha_f(r)+ K_1	 \int_{0}^{r} \frac{1}{s^{(\alpha-1)/2}}G_{z-y}^f(s)\mathrm{d}s.
\end{align}
It follows then from Gronwall's inequality that for all $y,z\in \R$ with $|y-z|<1$, $t>1$ and $r\in (0,T]$,
\begin{align}\label{step_49}
	& G_{z-y}^f(r)\leq \alpha_f(r)+	K_1\int_{0}^r\exp\left\{	K_1\int_{s}^r \frac{1}{q^{(\alpha-1)/2}} \mathrm{d}q \right\}\frac{\alpha_f(s)}{s^{(\alpha-1)/2}}\mathrm{d}s \nonumber\\
	& \lesssim \alpha_f(r)+   \int_{0}^r \frac{\alpha_f(s)}{s^{(\alpha-1)/2}}\mathrm{d}s .
\end{align}
Note that
\begin{align}\label{step_50}
	& \alpha_f(r)\lesssim r^{(\alpha-1)/2} \left( \frac{1}{\sqrt{tr}} \land 1 + \frac{1}{t}+\left(\frac{\epsilon_{tr}(h)}{\sqrt{r}}+ \frac{\sqrt{|y-z|}}{r^{3/4}} \right)\land \sqrt{t} +  |y-z|\right)
	 \nonumber\\
	 &\leq  r^{(\alpha-1)/2} \left( \frac{1}{t^{1/4}\sqrt{r}} + \frac{\sqrt{T}}{t^{1/4}\sqrt{r}}+\frac{\epsilon_{tr}(h)}{\sqrt{r}}+ \frac{\sqrt{|y-z|}}{r^{3/4}} + \frac{T^{3/4}\sqrt{ |y-z|}}{r^{3/4}}\right) \nonumber\\
	& \lesssim r^{(\alpha-1)/2} \left(\frac{1}{t^{1/4}\sqrt{r}}+ \frac{\epsilon_{tr}(h)}{\sqrt{r}}+ \frac{\sqrt{|y-z|}}{r^{3/4}} \right)
\end{align}
and
\begin{align}\label{step_51}
	&\int_0^r \frac{\epsilon_{ts}(h)}{\sqrt{s}}\mathrm{d}s \leq \int_0^r \frac{\widetilde{\epsilon}_{ts}(h)}{\sqrt{s}}\mathrm{d}s \leq \sup_{t>0} \epsilon_{t}(h)\int_0^{r/\sqrt{t}}\frac{1}{\sqrt{s}}\mathrm{d}s + \widetilde{\epsilon}_{\sqrt{t}r}(h)\int_0^T\frac{1}{\sqrt{s}}\mathrm{d}s \nonumber\\
	& \lesssim  \sqrt{T}\frac{1}{t^{1/4}} +\sqrt{T}\widetilde{\epsilon}_{\sqrt{t}r}(h)\lesssim \frac{1}{t^{1/4}}+ \widetilde{\epsilon}_{\sqrt{t}r}(h).
\end{align}
Therefore,
combining
\eqref{step_49}, \eqref{step_50}, \eqref{step_51} and the fact that $\epsilon_{tr}(h)\leq \widetilde{\epsilon}_{\sqrt{t}r}(h)$, we see that  for all $y,z\in \R$ with $|y-z|<1$, $t>1$ and $r\in (0,T]$,
\begin{align}\label{step_5}
	& G_{z-y}^f(r)\nonumber\\
	&\lesssim  r^{(\alpha-1)/2} \left(\frac{1}{t^{1/4}\sqrt{r}}+ \frac{\epsilon_{tr}(h)}{\sqrt{r}}+ \frac{\sqrt{|y-z|}}{r^{3/4}} \right)+ \left( \int_{0}^r\frac{1}{t^{1/4} \sqrt{s}} \mathrm{d}s+  \int_{0}^r \frac{\epsilon_{ts}(h)}{\sqrt{s}}\mathrm{d}s+
		\int_{0}^r
	 \frac{\sqrt{|y-z|}}{s^{3/4}}\mathrm{d}s \right) \nonumber\\
	&\lesssim r^{(\alpha-1)/2}\cdot \frac{1}{r^{3/4}} \left(\frac{T^{1/4}}{t^{1/4}}+ T^{1/4}\widetilde{\epsilon}_{\sqrt{t}r}(h)+ \sqrt{|y-z|} \right)\nonumber\\
	&\quad + r^{(\alpha-1)/2} \cdot \frac{T^{(5-2\alpha)/4}C(T)}{r^{3/4}}\left( \frac{2\sqrt{T}}{t^{1/4}} +  \left(\frac{1}{t^{1/4}}+ \widetilde{\epsilon}_{\sqrt{t}r}(h)\right)+4T^{1/4}\sqrt{|y-z|}\right)\nonumber\\
	&\lesssim r^{(\alpha-1)/2}\cdot \frac{1}{r^{3/4}}\left(\frac{1}{t^{1/4}}+\widetilde{\epsilon}_{\sqrt{t}r}(h)+\sqrt{|y-z|}\right),
\end{align}
which implies \eqref{step_34}.

Now we prove (ii).  Define
\[
\widehat{F}_q^f(r):= \sup_{x\in \R} \left| v_{g, h}^{(t)}(r,x)-v_{g, h}^{(t)}(r+q,x) \right|.
\]
Then combining
\eqref{Lemma:inequality}, Lemmas  \ref{lemma3}--\ref{Scaled-evolution-equation},
and $\frac{1}{\alpha-1}\geq 1$, we get that for all $y\in \R, q\in (0,1)$, $t>1$ and $r\in (0, T]$,
\begin{align}\label{step_7}
	& \left| v_{g, h}^{(t)}	(r,y)-v_{g, h}^{(t)}
	(r+q,y) \right| \leq  \frac{t^{\frac{1}{\alpha-1}}}{t^{\frac{2}{\alpha-1}-1}} \left(\mathbf{E}_y\left( h^2 (\sqrt{t}\xi_r^{(t)})\right) +\mathbf{E}_y\left( h^2 (\sqrt{t}\xi_{r+q}^{(t)})\right)  \right)\nonumber\\
	& \quad + \frac{1}{t^{\frac{1}{\alpha-1}}} \left(\mathbf{E}_y\left( g^2 (\xi_r^{(t)})\right) +\mathbf{E}_y\left( g^2 (\xi_{r+q}^{(t)})\right)  \right)+  \sqrt{t}\sup_{x\in\R} \left|\mathbf{E}_x(h(\sqrt{t}\xi_r^{(t)}))- \mathbf{E}_{x}(h(\sqrt{t}\xi_{r+q}^{(t)}))\right| \nonumber\\
	& \quad+  \sup_{x\in\R} \left|\mathbf{E}_x(g(\xi_r^{(t)}))- \mathbf{E}_{x}(g(\xi_{r+q}^{(t)}))\right| + \mathbf{E}_{y}\left(\int_r^{r+q} \psi^{(t)}\left(
	v_{g, h}^{(t)}
	(r+q-s,\xi_s^{(t)})\right)\mathrm{d}s\right) \nonumber\\
	&\quad + \left|\mathbf{E}_{y}\left(\int_0^r \psi^{(t)}\left(	v_{g, h}^{(t)}
	(r-s,\xi_s^{(t)})\right)\mathrm{d}s\right) -\mathbf{E}_{y}\left(\int_0^r \psi^{(t)}\left(	v_{g, h}^{(t)}
	(r+q-s,\xi_s^{(t)})\right)\mathrm{d}s\right)\right| \nonumber\\
	& \leq \mathbf{E}_y\left( h^2 (\sqrt{t}\xi_r^{(t)}) +\mathbf{E}_y\left( h^2 (\sqrt{t}\xi_{r+q}^{(t)})\right)  \right) + \frac{2}{t}\Vert g^2\Vert_\infty+ \sqrt{t}\sup_{x\in\R} \left|\mathbf{E}_x(h(\sqrt{t}\xi_r^{(t)}))- \mathbf{E}_{x}(h(\sqrt{t}\xi_{r+q}^{(t)}))\right|  \nonumber\\
	&\quad +  \sup_{x\in\R} \left|\mathbf{E}_x(g(\xi_r^{(t)}))- \mathbf{E}_{x}(g(\xi_{r+q}^{(t)}))\right| + \mathbf{E}_{y}\left(\int_r^{r+q} \psi^{(t)}\left(	v_{g, h}^{(t)}
	(r+q-s,\xi_s^{(t)})\right)\mathrm{d}s\right)\nonumber\\
	&\quad + C_\psi  \int_0^r \left(  \mathbf{E}_y\left( \left(	v_{g, h}^{(t)}
	(r-s,\xi_s^{(t)}) \right)^{\alpha-1}\right)+\mathbf{E}_y\left( \left(	v_{g, h}^{(t)}
	(r+q-s,\xi_s^{(t)}) \right)^{\alpha-1}\right) \right) \widehat{F}_q^f(r-s)\mathrm{d}s .
\end{align}
Combining Lemma \ref{lemma:small-effort}, \eqref{uniformly-bounded-2} and Lemma \ref{lemma4},
 for any $q\in (0,1), t>1$ and $r\in (0,T]$,
\begin{align}
	 &\widehat{F}_q^f(r)\lesssim \frac{1 }{\sqrt{tr}} \land 1+ \frac{1}{t}+ G_h(r,r+q;t) +\sqrt{q} \nonumber\\
	 & \quad + \frac{1}{\sqrt{r+q}}\int_r^{r+q} \frac{1}{(r+q-s)^{(\alpha-1)/2}}\mathrm{d}s + \frac{1}{r^{(\alpha-1)/2}} \int_0^r \widehat{F}_q^f(s)\mathrm{d}s\nonumber\\
	 &\lesssim \frac{1 }{\sqrt{tr}} \land 1+ \frac{1}{t}+ G_h(r,r+q;t) +\sqrt{q} + \frac{1}{\sqrt{r}}q^{(3-\alpha)/2}+ \frac{1}{r^{(\alpha-1)/2}} \int_0^r \widehat{F}_q^f(s)\mathrm{d}s.
\end{align}
Define $\widehat{G}_q^f(r):= r^{(\alpha-1)/2}\widehat{F}_q^f(r)$. Then  there exists a constant $K_2=K_2(g, h, T, \psi)$
such that for all $q\in (0,1), t>1$ and $r\in (0,T]$,
\begin{align}
	\widehat{G}_q^f(r)&\leq
	K_2  r^{(\alpha-1)/2}
		\left(\frac{1 }{\sqrt{tr}} \land 1+ \frac{1}{t}+ G_h(r,r+q;t) +\sqrt{q} + \frac{1}{\sqrt{r}}q^{(3-\alpha)/2}+ \frac{1}{r^{(\alpha-1)/2}} \int_0^r \widehat{F}_q^f(s)\mathrm{d}s\right) \nonumber\\
	&=:\widehat{\alpha}_f(r)+	K_2	\int_0^r\frac{ \widehat{G}_q^f(s)}{s^{(\alpha-1)/2}}\mathrm{d}s.
\end{align}
Therefore, by Gronwall's inequality, we get that for all $q\in (0,1), t>1$ and $r\in (0,T]$,
\begin{align}\label{step_64}
	 \widehat{G}_q^f(r)&\leq \widehat{\alpha}_f(r)+ K_2 \int_0^r \exp\left\{ 	 K_2 \int_s^r\frac{1}{q^{(\alpha-1)/2}}\mathrm{d}q \right\} \frac{\widehat{\alpha}_f(s)}{s^{(\alpha-1)/2}}\mathrm{d}s\nonumber\\
	 & \lesssim  \widehat{\alpha}_f(r)+  \int_0^r  \frac{\widehat{\alpha}_f(s)}{s^{(\alpha-1)/2}}\mathrm{d}s.
\end{align}
Using an argument similar to that
leading to
 \eqref{step_50} and  the fact that $(3-\alpha)/2 > 1/8$, we get that for all $q\in (0,1), t>1$ and $r\in (0,T]$,
\begin{align}\label{step_58}
	& \widehat{\alpha}_f(r) \lesssim r^{(\alpha-1)/2}  \left(  \frac{1 }{t^{1/4}\sqrt{r}} + \frac{\sqrt{T}}{t^{1/4}\sqrt{r}}+ G_h(r,r+q;t) +\frac{\sqrt{T}}{\sqrt{r}}q^{1/8} + \frac{1}{\sqrt{r}}q^{1/8} \right) \nonumber\\
	& \lesssim r^{(\alpha-1)/2}  \left(  \frac{1 }{t^{1/4}\sqrt{r}} +  G_h(r,r+q;t) + \frac{1}{\sqrt{r}}q^{1/8} \right).
\end{align}
 Moreover, by the definition of $G_h$, we know that for all $r\in (0, T], q\in(0,1)$ and $t>1$,
\begin{align}\label{step_10}
	& G_h(r,r+q;t)\leq \frac{2\widetilde{\epsilon}_{tr}(h)}{\sqrt{r}} +
	\ell(h)
	\left(\frac{q}{\sqrt{r(r+q)}(\sqrt{r}+\sqrt{r+q})} + \frac{1}{\sqrt{r+q}}\left(\sqrt{q}+\frac{\sqrt{q}}{r}\right)\right)\nonumber\\
	&\lesssim \frac{\widetilde{\epsilon}_{\sqrt{t}r}(h)}{\sqrt{r}} +\frac{q^{1/8}}{\sqrt{r^3}} + \frac{1}{\sqrt{r}}\left(1+\frac{1}{r}\right)q^{1/8}\lesssim \frac{1}{r^{3/2}}\left(\widetilde{\epsilon}_{\sqrt{t}r}(h)+ q^{1/8}\right).
\end{align}
Combining \eqref{step_58} and \eqref{step_10}, we see that  for all $t>1$, $r\in (0, T]$ and $q\in (0,1)$,
\begin{align}\label{step_59}
	\widehat{\alpha}_f(r)\lesssim r^{(\alpha-1)/2} \cdot \frac{1}{r^{3/2}}\left(\frac{1}{t^{1/4}}+ q^{1/8}+ \widetilde{\epsilon}_{\sqrt{t}r}(h)\right).
\end{align}
Using \eqref{step_58}, we get that for all $t>1$, $r\in (0, T]$ and $q\in (0,1)$,
\begin{align}\label{step_60}
	& \int_0^r  \frac{\widehat{\alpha}_f(s)}{s^{(\alpha-1)/2}}\mathrm{d}s \lesssim  \int_0^r   \frac{1}{t^{1/4}\sqrt{s}}\mathrm{d}s + \int_0^r  G_h(s,s+q;t)\mathrm{d}s+\int_0^r  \frac{q^{1/8}}{\sqrt{s}}\mathrm{d}s \nonumber\\
	& \lesssim    \frac{1}{t^{1/4}} + \int_0^r  G_h(s,s+q;t)\mathrm{d}s+  q^{1/8}.
\end{align}
Note that by \eqref{step_51},
\begin{align}\label{step_8}
	&\int_0^r  G_h(s,s+q;t)\mathrm{d}s \leq  \int_0^r \frac{2\widetilde{\epsilon}_{ts}(h)}{\sqrt{s}}\mathrm{d}s+ \int_0^r \left(\frac{1}{\sqrt{s}}-\frac{1}{\sqrt{s+q}} + \frac{1}{\sqrt{s+q}}\left(\sqrt{q}+1-e^{-\frac{\sqrt{q}}{s}}\right) \right)\mathrm{d}s\nonumber\\
	&\lesssim \left(\frac{1}{t^{1/4}}+\widetilde{\epsilon}_{\sqrt{t}r}(h)\right)+\int_0^r \left(\frac{1}{\sqrt{s}}-\frac{1}{\sqrt{s+q}} + \frac{1}{\sqrt{s+q}}\left(\sqrt{q}+1-e^{-\frac{\sqrt{q}}{s}}\right) \right)\mathrm{d}s
\end{align}
and  that
\begin{align}\label{step_61}
	&\int_0^r \left(\frac{1}{\sqrt{s}}-\frac{1}{\sqrt{s+q}} + \frac{1}{\sqrt{s+q}}\left(\sqrt{q}+1-e^{-\frac{\sqrt{q}}{s}}\right) \right)\mathrm{d}s \nonumber\\
	& = 2(\sqrt{r}-\sqrt{r+q}+\sqrt{q})+ \int_0^r \left( \frac{1}{\sqrt{s+q}}\left(\sqrt{q}+1-e^{-\frac{\sqrt{q}}{s}}\right) \right)\mathrm{d}s \nonumber\\
	&\leq 2\sqrt{q}+ 2\int_0^{q^{1/4}}\frac{1}{\sqrt{s+q}}\mathrm{d}s + \int_{q^{1/4}}^{r\vee q^{1/4} } \frac{1}{\sqrt{s+q}}(\sqrt{q}+ \frac{\sqrt{q}}{s})\mathrm{d}s\nonumber\\
	&\leq 2q^{1/8}+4\sqrt{q+q^{1/4}}+ (\sqrt{q}+q^{1/4})\int_0^T \frac{1}{\sqrt{s}}\mathrm{d}s\lesssim q^{1/8}.
\end{align}
Combining \eqref{step_60}, \eqref{step_8} and \eqref{step_61}, we get that for all $t>1$, $q\in (0,1)$ and  $r\in (0, T]$,
\begin{align}\label{step_63}
	\int_0^r  \frac{\widehat{\alpha}_f(s)}{s^{(\alpha-1)/2}}\mathrm{d}s  \lesssim  \frac{1}{t^{1/4}}+ \widetilde{\epsilon}_{\sqrt{t}r}(h)+ q^{1/8} \lesssim r^{(\alpha-1)/2} \times \frac{1}{r^{3/2}}\left(\frac{1}{t^{1/4}}+ q^{1/8}+ \widetilde{\epsilon}_{\sqrt{t}r}(h)\right).
\end{align}
Now combining \eqref{step_64}, \eqref{step_59} and \eqref{step_63},  we get \eqref{step_35}.
\hfill$\Box$
\bigskip

By \eqref{upper-bound-of-v-theta-A}, we see that for any $r>0$,
$h\in {\rm DRI}^+_c(\R)$
and $g\in B_{Lip}^+(\R)$, we have $\sup_{y\in \R, t>1} v_{g, h}^{(t)}(r,y)<\infty$.
Therefore, by a diagonalization argument, for any sequence of positive reals increasing to $\infty$, we can find a subsequence $\{t_k:k\in \N\}$ such that $\lim_{k\to\infty} t_k=\infty$ and that the following limit exists:
\begin{align}\label{subsequence-limit}
	\lim_{k\to\infty}
		v_{g, h}^{(t_k)}
	(r,y)=: v_{g,h}^X(r,y),\quad\mbox{ for all } r\in \Q_+:= (0,\infty)\cap \Q, y\in \Q.
\end{align}
Of course, the choice of $\{t_k\}$ may depend on the functions $h$ and $g$.
Using \eqref{limit0}, taking $t=t_k$ in Proposition \ref{prop1} first and then letting $k\to\infty$, we see that $v_{g,h}^X(r,y)$ is continuous in $\Q_+\times \Q$.
Now for each $r>0, y\in \R$,  for any sequence $\{(r_m, y_m),m\in \N\}\subset \Q_+\times \Q$ with $r_m \to r$ and $y_m\to y$ as $m\to\infty$, we see that the sequence of $\left\{v_{g,h}^X(r_m,y_m), m\in \N\right\}$ is a Cauchy sequence. Therefore, for each $r>0$ and $y\in \R$, we can define
\begin{align}
	v_{g,h}^X(r,y):= \lim_{(r_m, y_m)\in \Q_+\times \Q, (r_m, y_m)\to (r,y)} v_{g,h}^X(r_m,y_m).
\end{align}

Our next result shows that \eqref{subsequence-limit}   also holds for all $r>0, y\in \R$.

\begin{lemma}
	The limits \eqref{subsequence-limit}  holds for all $r>0$ and $y\in \R$.
\end{lemma}
\textbf{Proof: } Let $(r_m, y_m)\in \Q_+\times \Q$ be such that $(r_m, y_m)\to (r,y)$. Without loss of generality, we assume that $2r>r_m> \frac{1}{2}r$ for all $m$. Then by Proposition \ref{prop1} with $T=2r$,
 there exist positive constants $N_1=N_1(g,h,r,\psi)$ and $N_2=N_2(g,h,r,\psi)$ such that
\begin{align}
	& \left|  v_{g,h}^X(r,y) -	v_{g,h}^{(t_k)}(r,y) 	\right| \leq \left|  v_{g,h}^X(r,y) -v_{g,h}^X(r_m,y_m)\right| +\left| v_{g,h}^X(r_m,y_m) -	v_{g, h}^{(t_k)}	(r_m,y_m) \right| \nonumber\\
	& \quad + \left|	v_{g, h}^{(t_k)}	(r_m,y_m)-	v_{g, h}^{(t_k)}	(r,y)\right|\nonumber\\
	& \leq \left|  v_{g,h}^X(r,y) -v_{g,h}^X(r_m,y_m)\right| +\left| v_{g,h}^X(r_m,y_m) -
	v_{g, h}^{(t_k)}	(r_m,y_m) \right| \nonumber\\
	& \quad + 	\frac{N_1 2^{3/4}}{r^{3/4}}
	 \left(\widetilde{\epsilon}_{\sqrt{t}r/2}(h)+ \sqrt{|r_m-r|}+ \frac{1}{t_k^{1/4}}\right)+ 
	\frac{N_22^{3/2}}{r^{3/2}}
	\left(\widetilde{\epsilon}_{\sqrt{t}r/2}(h)+ |y_m-y|^{1/8}+\frac{1}{t_k^{1/4}}\right).
\end{align}
Letting $k\to\infty$ in the inequality above and using  \eqref{subsequence-limit}, we get
\begin{align}
	& \limsup_{k\to\infty} \left|  v_{g,h}^X(r,y) -	v_{g, h}^{(t_k)}	(r,y) \right|
	\nonumber\\ &\leq \left|  v_{g,h}^X(r,y) -v_{g,h}^X(r_m,y_m)\right| +  \frac{2^{3/4}}{r^{3/4}}
		N_1 \sqrt{|r_m-r|} + N_2\frac{2^{3/2}|y_m-y|^{1/8}}{r^{3/2}}.
\end{align}
Letting $m\to\infty$, we arrive at the desired result. This completes the proof of the lemma.

\hfill$\Box$

Define
\begin{align}\label{Def-of-eta-theta-A}
      \eta_{g,h}(\mathrm{d}x):=	\ell(h)\delta_0(\mathrm{d}x) + g(x)\mathrm{d}x.
\end{align}

\begin{prop}\label{prop2}
	Any subsequential limit $v_{g,h}^X(r,y)$ of $\{ v_{g, h}^{(t)}	(r,y)\}$ is equal to the solution
   $v_{(\emptyset, \eta_{g,h})}^X(r,y)$ of \eqref{PDE} with initial trace $(\emptyset, \eta_{g,h})$
whose probabilistic representation  is given in Proposition \ref{Probabilistic-representation}.
\end{prop}
\textbf{Proof: }
By the uniqueness of solutions to  \eqref{PDE} and Proposition \ref{Probabilistic-representation}, we only need to prove that any subsequential limit $v_{g,h}^X(r,y)$ is the solution of \eqref{PDE}.
We divide the proof to two steps.

\textbf{Step 1:} In this step we
derive the integral  equation for any subsequential limit $v_{g,h}^X(r,y)$.
Noticing that $\frac{1}{\alpha-1}\geq 1$,
by \eqref{Lemma:inequality},
for any $r>0$ and $y\in \R$, it holds that
\begin{align}\label{step_69}
	& \limsup_{t\to\infty} \left|t^{\frac{1}{\alpha-1}}\mathbf{E}_{y}\left(1-\exp\left\{-\frac{1}{t^{\frac{1}{\alpha-1}-\frac{1}{2}}}h(\sqrt{t}\xi_{r}^{(t)}) -\frac{1}{t^{\frac{1}{\alpha-1}}}g(\xi_r^{(t)})\right\}\right)  - \sqrt{t}\mathbf{E}_y\left(h(\sqrt{t}\xi_{r}^{(t)}) \right)- \mathbf{E}_y\left(g(\xi_r^{(t)})\right)\right|\nonumber\\
	&\leq \limsup_{t\to\infty } \left(\frac{t^{\frac{1}{\alpha-1}}}{t^{\frac{2}{\alpha-1}-1}} \mathbf{E}_y\left(h^2(\sqrt{t}\xi_r^{(t)})\right) + \frac{t^{\frac{1}{\alpha-1}}}{t^{\frac{2}{\alpha-1}}} \mathbf{E}_y\left(g^2(\xi_r^{(t)})\right)\right)\nonumber\\
	&\leq  \limsup_{t\to\infty } \left( \Vert h\Vert_\infty \times \frac{C(g,h)}{\sqrt{tr}}+ \frac{1}{t} \Vert g^2\Vert_\infty\right) =0,
\end{align}
where $C(g,h)$ is defined in \eqref{uniformly-bounded-2}.
Note that by \eqref{scaling},
\[
\sqrt{t}\mathbf{E}_{y}\left(h(\sqrt{t}\xi^{(t)}_r)\right)-\frac{\ell(h)}{\sqrt{r}}\phi\left(\frac{y}{\sqrt{r}}\right)
=\frac1{\sqrt{r}}\left(\sqrt{tr}\mathbf{E}_{\sqrt{t}y}h(\xi_{tr})-\ell(h)\phi\left(\frac{\sqrt{t}y}{\sqrt{tr}}\right)\right).
\]
Thus by Lemma \ref{lemma1},
\begin{align}\label{e:rs}
\lim_{t\to\infty}\sqrt{t}\mathbf{E}_{y}\left(h(\sqrt{t}\xi^{(t)}_r)\right)=\frac{\ell(h)}{\sqrt{r}}\phi\left(\frac{y}{\sqrt{r}}\right).
\end{align}
By the central limit theorem we know that $\xi_r^{(t)}$ converges weakly to $B_r$. Then, combining \eqref{step_69} and \eqref{e:rs},
we get that
\begin{align}\label{step_68}
	& \lim_{t\to\infty} t^{\frac{1}{\alpha-1}}\mathbf{E}_{y}\left(1-\exp\left\{-\frac{1}{t^{\frac{1}{\alpha-1}-\frac{1}{2}}}h(\sqrt{t}\xi_{r}^{(t)}) -\frac{1}{t^{\frac{1}{\alpha-1}}}g(\xi_r^{(t)})\right\}\right)  \nonumber\\
	& =
	\frac{ \ell(h)}
	{\sqrt{r}}\phi\left(\frac{y}{\sqrt{r}}\right)+\mathbf{E}_y\left(g(B_r)\right).
\end{align}
Now letting $t=t_k$ in Lemma \ref{Scaled-evolution-equation}  first and then $k\to\infty$, by \eqref{step_68},   for any $r>0$ and $y\in \R$,
\begin{align}\label{step_12}
	v_{g,h}^X(r,y) = \frac{
		 \ell(h)}
	{\sqrt{r}}\phi\left(\frac{y}{\sqrt{r}}\right)+\mathbf{E}_y\left(g(B_r)\right)-\lim_{k\to\infty} \mathbf{E}_{y}\left(\int_0^r \psi^{(t_k)}\left(	v_{g, h}^{(t_k)}	(r-s,\xi_s^{(t_k)})\right)\mathrm{d}s\right).
\end{align}
Combining  \eqref{Properties-of-phi} and \eqref{upper-bound-of-v-theta-A}, for any $\varepsilon>0$, there exists $N\in \N$ such that for all $s\in (0, r-\varepsilon)$ and $y\in \R$,
\begin{align}\label{step_14}
	(1-\varepsilon) \varphi \left(	v_{g, h}^{(t_k)}	(r-s, y)\right)\leq \psi^{(t_k)}\left(	v_{g, h}^{(t_k)}
	(r-s,y)\right) \leq (1+\varepsilon) \varphi \left(	v_{g, h}^{(t_k)}	(r-s, y)\right),
\end{align}
where $\varphi$ is the function defined in \eqref{Stable-Branching-mechanism}.
Thus, for $k>N$, we have
\begin{align}\label{Lower-bound-1}
	& \mathbf{E}_{y}\left(\int_0^r \psi^{(t_k)}\left(	v_{g, h}^{(t_k)}
	(r-s,\xi_s^{(t_k)})\right)\mathrm{d}s\right) \geq \mathbf{E}_{y}\left(\int_0^{r-\varepsilon} \psi^{(t_k)}\left(
	v_{g, h}^{(t_k)}(r-s,\xi_s^{(t_k)})\right)\mathrm{d}s\right)\nonumber\\
	& \geq  (1-\varepsilon)\mathbf{E}_{y}\left(\int_0^{r-\varepsilon} \varphi\left(	v_{g, h}^{(t_k)}	(r-s,\xi_s^{(t_k)})\right)\mathrm{d}s\right).
\end{align}
Similarly, combining  Lemma \ref{lemma3}, \eqref{upper-bound-of-v-theta-A} and \eqref{step_14}, we get that for $k>N$,
\begin{align}\label{Upper-bound-1}
	& \mathbf{E}_{y}\left(\int_0^r \psi^{(t_k)}\left(		 v_{g, h}^{(t_k)}
	(r-s,\xi_s^{(t_k)})\right)\mathrm{d}s\right) \nonumber\\
	& \leq  (1+\varepsilon)\mathbf{E}_{y}\left(\int_0^{r-\varepsilon} \varphi\left(	v_{g, h}^{(t_k)}
	(r-s,\xi_s^{(t_k)})\right)\mathrm{d}s\right) +C_\psi \int_{r-\varepsilon}^r \left(\frac{	C}{\sqrt{r-s}}\right)^\alpha\mathrm{d}s \nonumber\\
	& = (1+\varepsilon)\mathbf{E}_{y}\left(\int_0^{r-\varepsilon} \varphi\left(	v_{g, h}^{(t_k)}
	(r-s,\xi_s^{(t_k)})\right)\mathrm{d}s\right) + f(\varepsilon),
\end{align}
where $C$ is a positive constant, and $f(\varepsilon)$ is a function of  $\varepsilon$ satisfying  $\lim_{\varepsilon\to 0} f(\varepsilon)=0$.
We claim that for any $r>\varepsilon>0$ and $y\in \R$,
\begin{align}\label{Claim-1}
	\lim_{k\to\infty} \mathbf{E}_{y}\left(\int_0^{r-\varepsilon} \varphi\left(
		v_{g, h}^{(t_k)}
	(r-s,\xi_s^{(t_k)})\right)\mathrm{d}s\right) = \mathbf{E}_{y}\left(\int_0^{r-\varepsilon} \varphi\left(     v_{g, h}^X	(r-s,B_s)\right)\mathrm{d}s\right).
\end{align}

To prove \eqref{Claim-1}, fix two large constants $R$ and
$T>r$,
since $v_{g,h}^X(r,y)$ is continuous in $(r,y)\in [\varepsilon, T]\times [-R, R]$, for any
$\gamma_0\in (0,1)$,
 there exist $L\in\N$ and $r_0=\varepsilon<...< r_L=T$, $y_0:=-R<...<y_L=R$ such that
 $\max_{i\in \{1,...,L\}} |r_i-r_{i-1}|<\gamma_0, \max_{i\in \{1,...,L\}} |y_i-y_{i-1}|<\gamma_0$ and that
\begin{align}\label{Proof-Claim1-1}
	\max_{i, j\in \{1,...,L\}} \max_{r\in [r_{i-1}, r_i], y\in [y_{j-1}, y_j]} \left| v_{g,h}^X(r_i,y_j) - v_{g,h}^X(r,y) \right| <
	\gamma_0.
\end{align}
Now we take $N_*$ sufficiently large so that when $k>N_*$,
\begin{align}\label{Proof-Claim1-2}
	\max_{r\in \{r_0,..., r_L\}, y\in \{y_0,...,y_L\}} \left| v_{g, h}^{(t_k)}
	(r,y) -v_{g,h}^X(r,y)  \right|<
	\gamma_0.
\end{align}
Combining Proposition \ref{prop1}, \eqref{Proof-Claim1-1} and \eqref{Proof-Claim1-2},
we get that, if  $k>N_*$, then for any $r\in [\varepsilon, T], y\in [-R, R]$, suppose that $r_{i_0-1}\leq r\leq r_{i_0}$ and $y_{j_0-1}\leq y\leq y_{j_0}$ for some $i_0, j_0\in \{1,...,L\}$,
\begin{align}
	& \left|	v_{g, h}^{(t_k)}(r,y) -v_{g,h}^X(r,y) \right|<
	\gamma_0
	+ \left|v_{g, h}^{(t_k)} (r,y) -v_{g,h}^X(r_{i_0},y_{j_0}) \right|\nonumber\\
	& <
	2\gamma_0
	+ \left|	v_{g, h}^{(t_k)}	(r,y) -	v_{g,h }^{(t_k)}	(r_{i_0},y_{j_0}) \right| \nonumber\\
	& \leq 2 \gamma_0 +
	\frac{N_1}{\varepsilon^{3/4}}\left(\widetilde{\epsilon}_{\sqrt{t_k}\varepsilon}(h)+\sqrt{|y-y_{j_0}|}+ \frac{1}{t_k^{1/4}}\right) +\frac{N_2}{\varepsilon^{3/2}} \left(\widetilde{\epsilon}_{\sqrt{t_k}\varepsilon}(h)+ \left|r-r_{j_0}\right|^{1/8}+ \frac{1}{t_k^{1/4}}\right)\nonumber\\
	& \leq 2 \gamma_0^{1/8} +
  \left(\frac{N_1}{\varepsilon^{3/4}}+ \frac{N_2} {\varepsilon^{3/2}}\right)  \left(\widetilde{\epsilon}_{\sqrt{t_k}\varepsilon}(h)+
  \gamma_0^{1/8}
  +\frac{1}{t_k^{1/4}}\right) .
\end{align}
Note that $\lim_{t\to\infty}\widetilde{\epsilon}_{\sqrt{t}\varepsilon}(h)=0$ by \eqref{limit0}. Therefore,
for any
$\gamma_0\in (0,1)$,
 there exists a constant $C'=C'(\varepsilon, T, R)$ such that when $k$ is large enough,
\begin{align}
	& \left|	v_{g, h}^{(t_k)}	(r,y) -v_{g,h}^X(r,y) \right| \leq C'
	\gamma_0^{1/8},
\end{align}
for any $r\in [\varepsilon, T]$ and $y\in [-R, R]$. Therefore, combining \eqref{upper-bound-of-v-theta-A} and the fact that $|\varphi(u)-\varphi(v)|\leq C_{\varphi, \varepsilon} |u-v|$ for all
$u,v \in [0, \frac{C}{\sqrt{\varepsilon}}]$,
we obtain that for any
$\gamma\in (0,1) $,
 $R>1$ and $T>r$,
\begin{align}
	&\limsup_{k\to\infty}\left| \mathbf{E}_{y}\left(\int_0^{r-\varepsilon} \varphi\left(
		v_{g, h}^{(t_k)}
	(r-s,\xi_s^{(t_k)})\right)\mathrm{d}s\right) -  \mathbf{E}_{y}\left(\int_0^{r-\varepsilon} \varphi\left(v_{g,h}^X(r-s,\xi_s^{(t_k)})\right)\mathrm{d}s\right) \right|\nonumber\\
	& \leq 2 \limsup_{k\to\infty}\int_0^{r-\varepsilon} \varphi\left(\frac{
		C}{\sqrt{r-s}}\right) \mathbf{P}_y\left(|\xi_s^{(t_k)}|>R\right)\mathrm{d}s\nonumber\\
	&\quad\quad +C_{\varphi,\varepsilon} \limsup_{k\to\infty} \mathbf{E}_{y}\left(\int_0^{r-\varepsilon} \left|
		v_{g, h}^{(t_k)}
	(r-s,\xi_s^{(t_k)}) -v_{g,h}^X(r-s,\xi_s^{(t_k)})\right|1_{\{\xi_s^{(t_k)}\in [-R,R]\}}\mathrm{d}s\right)\nonumber\\
	& \leq 2 \varphi\left(\frac{
		C	}{\sqrt{\varepsilon}}\right) \int_0^{r-\varepsilon}\frac{s+y^2}{R^2}\mathrm{d}s + C_{\varphi,\varepsilon} C'
		\gamma_0^{1/8} (r-\varepsilon) \stackrel{\gamma_0\downarrow 0, R\uparrow \infty}{\longrightarrow}0.
\end{align}
By the functional central limit theorem we know that
$(\xi_s^{(t)})_{0<s\leq r-\varepsilon}$ converges weakly to $(B_s)_{0< s\leq r-\varepsilon}$,
thus \eqref{Claim-1} is valid.

Plugging \eqref{Claim-1} into \eqref{Lower-bound-1} and \eqref{Upper-bound-1}, we conclude that for any $\varepsilon>0$,
\begin{align}
	& (1-\varepsilon)\mathbf{E}_{y}\left(\int_0^{r-\varepsilon} \varphi\left(v_{g,h}^X(r-s,B_s)\right)\mathrm{d}s\right)\nonumber\\
	& \leq \liminf_{k\to\infty} \mathbf{E}_{y}\left(\int_0^r \psi^{(t_k)}\left(	v_{g, h}^{(t_k)}
	(r-s,\xi_s^{(t_k)})\right)\mathrm{d}s\right) \leq \limsup_{k\to\infty} \mathbf{E}_{y}\left(\int_0^r \psi^{(t_k)}\left(	v_{g, h}^{(t_k)}	(r-s,\xi_s^{(t_k)})\right)\mathrm{d}s\right) \nonumber\\
	& \leq  (1+\varepsilon)\mathbf{E}_{y}\left(\int_0^{r-\varepsilon} \varphi\left(v_{g,h}^X(r-s,B_s)\right)\mathrm{d}s\right)+ f(\varepsilon).
\end{align}
Letting $\varepsilon\downarrow0$ in the inequality above and plugging the resulting inequality into \eqref{step_12}, we get that for any $r>0$ and $y\in\R$,
\begin{align}\label{Evolution-equation}
	v_{g,h}^X(r,y) = \frac{	\ell(h)	}{\sqrt{r}}\phi\left(\frac{y}{\sqrt{r}}\right)+\mathbf{E}_y(g(B_r))-\mathbf{E}_{y}\left(\int_0^{r} \varphi\left(v_{g,h}^X(r-s,B_s)\right)\mathrm{d}s\right).
\end{align}

\textbf{Step 2:} In this step we show that  \eqref{Evolution-equation} is equivalent to \eqref{PDE} with initial trace $(\emptyset, \eta_{g,h})$. Combining \eqref{Evolution-equation} and the Markov property, we see that for $w\in (0,r)$,
\begin{align}\label{step_15}
		& v_{g,h}^X(r,y) + \mathbf{E}_{y}\left(\int_0^{w} \varphi\left(v_{g,h}^X(r-s,B_s)\right)\mathrm{d}s\right)\nonumber\\
		& =  \frac{	\ell(h)	}{\sqrt{r}}\phi\left(\frac{y}{\sqrt{r}}\right)+\mathbf{E}_y(g(B_r))- \mathbf{E}_{y}\left(\int_w^{r} \varphi\left(v_{g,h}^X(r-s,B_s)\right)\mathrm{d}s\right)\nonumber\\
		& = \frac{\ell(h)	}{\sqrt{r}}\phi\left(\frac{y}{\sqrt{r}}\right) + \mathbf{E}_{y}\left( \mathbf{E}_{B_w}(g(B_{r-w}))- \mathbf{E}_{B_w}\left( \int_0^{r-w} \varphi\left(v_{g,h}^X(r-s,B_s)\right)\mathrm{d}s\right)\right)\nonumber\\
		& =  \frac{	\ell(h)	}{\sqrt{r}}\phi\left(\frac{y}{\sqrt{r}}\right) - \mathbf{E}_{y}\left(\frac{	\ell(h)}
		{\sqrt{r-w}}\phi\left(\frac{B_w}{\sqrt{r-w}}\right) \right) +\mathbf{E}_{y}\left( v_{g,h}^X(r-w, B_w) \right).
\end{align}
Routine computations yield that
\begin{align}\label{step_16}
	& \frac{\ell(h)}{\sqrt{r}}\phi\left(\frac{y}{\sqrt{r}}\right) - \mathbf{E}_{y}\left(\frac{
		\ell(h)	}{\sqrt{r-w}}\phi\left(\frac{B_w}{\sqrt{r-w}}\right) \right)\nonumber\\
	& =  \frac{	\ell(h)}{\sqrt{2\pi}} \left( \frac{1}{\sqrt{r}}e^{-y^2/(2r)}- \frac{1}{\sqrt{2\pi w(r-w)}} \int \exp\left\{ - \frac{rz^2}{2w(r-w)}+\frac{zy}{w}-\frac{y^2}{2w}\right\} \mathrm{d}z\right)\nonumber\\
	& =  \frac{	\ell(h)	}{\sqrt{2\pi}} \left( \frac{1}{\sqrt{r}}e^{-y^2/(2r)}- \frac{1}{\sqrt{2\pi w(r-w)}} \int \exp\left\{ - \frac{r}{2w(r-w)}\left(z-\frac{r-w}{r}y\right)^2-\frac{y^2}{2r}\right\} \mathrm{d}z\right)\nonumber\\
	& =0.
\end{align}
Therefore, combining \eqref{step_15} and  \eqref{step_16}, we conclude that
\begin{align}\label{step_17}
	 v_{g,h}^X(r,y) + \mathbf{E}_{y}\left(\int_0^{w} \varphi\left(v_{g,h}^X(r-s,B_s)\right)\mathrm{d}s\right) = \mathbf{E}_{y}\left( v_{g,h}^X(r-w, B_w) \right).
\end{align}
For any fixed $w>0$, set $u(r, y):= v_{g,h}^X(r+w,y)$, then we see from \eqref{step_17} that $u$ solves \eqref{Evolution-cumulant-semigroup} with $f= v_{g,h}^X(w,\cdot)$.  Now it suffices to check the boundary condition. By \eqref{Evolution-equation},  for any $j\in C_c^+(\R)$ and any $r>0$,
\begin{align}
	&\int j(y)	v_{g,h}^X(r,y)\mathrm{d}y\nonumber\\
	& =	\ell(h) \int j(y)\frac{1}{\sqrt{r}}\phi\left(\frac{y}{\sqrt{r}}\right)+ \int j(y)\mathbf{E}_y\left(g(B_r)\right)\mathrm{d}s- \int j(y)\mathbf{E}_{y}\left(\int_0^{r} \varphi\left(v_{g,h}^X(r-s,B_s)\right)\mathrm{d}s\right)\mathrm{d}y\nonumber\\
	& =	\ell(h)	\mathbf{E}_0\left(j(B_r)\right)+ \int j(y)\mathbf{E}_y\left(g(B_r)\right)\mathrm{d}s- \int j(y)\mathbf{E}_{y}\left(\int_0^{r} \varphi\left(v_{g,h}^X(r-s,B_s)\right)\mathrm{d}s\right)\mathrm{d}y.
\end{align}
Since $\lim_{r\downarrow 0}\mathbf{E}_0\left(j(B_r)\right) = j(0)$ and $\lim_{r\downarrow 0}  \int j(y)\mathbf{E}_y\left(g(B_r)\right)\mathrm{d}s= \int j(y)g(y)\mathrm{d}y$ by the dominated convergence theorem,
to prove the desired result, we only need to prove that
\begin{align}\label{step_18}
	\lim_{r\downarrow 0} \int j(y)\mathbf{E}_{y}\left(\int_0^{r} \varphi\left(v_{g,h}^X(r-s,B_s)\right)\mathrm{d}s\right)\mathrm{d}y=0.
\end{align}
Combining \eqref{upper-bound-of-v-theta-A} and the definition of $\varphi$, we see that
\begin{align}\label{step_19}
	& \int j(y)\mathbf{E}_{y}\left(\int_0^{r} \varphi\left(v_{g,h}^X(r-s,B_s)\right)\mathrm{d}s\right)\mathrm{d}y\nonumber\\
	&\lesssim  \int j(y)\mathbf{E}_{y}\left(\int_0^{r} \frac{1}{\left(\sqrt{r-s}\right)^{\alpha-1}}v_{g,h}^X(r-s,B_s)\mathrm{d}s\right)\mathrm{d}y\nonumber\\
	& = \int_0^{r} \frac{1}{\left(r-s\right)^{(\alpha-1)/2}} \int j(y) \mathbf{E}_{y}\left(v_{g,h}^X(r-s,B_s)\right)  \mathrm{d}y  \mathrm{d}s\nonumber\\
	& \leq  \int_0^{r} \frac{1}{\left(r-s\right)^{(\alpha-1)/2}}  \int j(y) \mathbf{E}_{y}\left(\frac{
	\ell(h)	}{\sqrt{r-s}}\phi\left(\frac{B_s}{\sqrt{r-s}}\right)+ \mathbf{E}_{B_s}(g(B_{r-s}))\right)  \mathrm{d}y  \mathrm{d}s,
\end{align}
where in  the last inequality we used \eqref{Evolution-equation}.
Combining  \eqref{step_16} and \eqref{step_19}, we get
\begin{align}
	& \int j(y)\mathbf{E}_{y}\left(\int_0^{r} \varphi\left(v_{g,h}^X(r-s,B_s)\right)\mathrm{d}s\right)\mathrm{d}y\nonumber\\
	&\lesssim \int_0^{r} \frac{1}{\left(r-s\right)^{(\alpha-1)/2}} \int j(y)\left(\frac{	\ell(h)
	}{\sqrt{r}}\phi\left(\frac{y}{\sqrt{r}}\right)  +\Vert g\Vert_\infty  \right) \mathrm{d}y\mathrm{d}s\nonumber\\
	& =   \int_0^{r} \frac{1}{s^{(\alpha-1)/2}}  \mathrm{d}s \left(	\ell(h)
	\mathbf{E}_0(j(B_r))+ \Vert g\Vert_\infty\int j(y)\mathrm{d}y\right) \nonumber\\
	&\lesssim \left(	\ell(h)	\Vert j\Vert_\infty +	\ell(j)	\Vert g\Vert_\infty \right)\int_0^{r} \frac{1}{s^{(\alpha-1)/2}}  \mathrm{d}s\stackrel{r\downarrow 0}{\longrightarrow} 0,
\end{align}
which implies the desired result.
\hfill$\Box$

\textbf{Proof of Theorem \ref{thm1}:} Theorem \ref{thm1} follows directly from Proposition \ref{prop2}.

\hfill$\Box$

Now we are going to prove Theorem \ref{thm2}. Before the proof, we need to prove an upper bound for maximal position $M_{t}:=\max_{u\in N(t)}X_u(t)$ and the minimal position  $M_t^-:=\min_{u\in N(t)}X_u(t)$ of all the particle alive at time $t$ with the convention that $M_t=-\infty$ and $M_t^-=\infty$ when $N(t)=\emptyset.$
In the next lemma, we will need
{\bf(H4')}
 to control the overshoot of L\'{e}vy process.

\begin{lemma}\label{Technical-Lemma}
		Assume {\bf(H1)}, {\bf(H2)}, {\bf(H3)} and {\bf(H4')} hold. For
	any $q,\delta>0$, there exist constants $C(q), T(\delta)\in (0,\infty)$ such that for $t>T(\delta)$,
	\begin{align}
		 t^{\frac{1}{\alpha-1}}\P_0\left( M_{t\delta}> q\sqrt{t}\right)\leq C(q)\delta \quad\mbox{and}\quad  t^{\frac{1}{\alpha-1}}\P_0\left( M_{t\delta}^- <- q\sqrt{t}\right)\leq C(q)\delta.
	\end{align}
\end{lemma}
\textbf{Proof: }  We only prove the first inequality, the proof of the second one is similar. Set
\[
Q^{(t)}(r,x):= t^{\frac{1}{\alpha-1}}\P_{\sqrt{t}x}\left( M_{tr}> 0\right)= \lim_{\theta\uparrow\infty} t^{\frac{1}{\alpha-1}}\left(1-\E_{\sqrt{t}x}\left(\exp\left\{-\theta Z_{tr}((0,\infty)) \right\}\right)\right).
\]
Then $t^{\frac{1}{\alpha-1}}\P_0\left( M_{t\delta}> q\sqrt{t}\right)= Q^{(t)}(\delta, -q)$, and we only need to prove that there exist constants $C(q), T(\delta)\in (0,\infty)$ such that for $t>T(\delta)$,
\begin{equation}\label{toprove}
Q^{(t)}(\delta, -q)\leq C(q)\delta.
\end{equation}
By Lemma \ref{Scaled-evolution-equation} (with $h=\theta 1_{(0,\infty)}, g=0$ first and then $\theta\uparrow\infty$), we see that $Q^{(t)}(r,x)$ solves
\begin{align}
	Q^{(t)}(r,x)=t^{\frac{1}{\alpha-1}}\mathbf{P}_{x}\left(\xi_{r}^{(t)}>0\right) - \mathbf{E}_{x}\left(\int_0^r \psi^{(t)}\left(Q^{(t)}(r-s,\xi_s^{(t)})\right)\mathrm{d}s\right).
\end{align}
By the Markov property, for any $w<r$, the above equation can also be rewritten by
\begin{align}
	Q^{(t)}(r,x)=\mathbf{E}_x\left(Q^{(t)}(r-w, \xi_{w}^{(t)})\right)- \mathbf{E}_{x}\left(\int_0^w \psi^{(t)}\left(Q^{(t)}(r-s,\xi_s^{(t)})\right)\mathrm{d}s\right).
\end{align}
It follows from the Feynman-Kac formula that for any $0<w<r$,
\begin{align}\label{Feynman-Kac}
	Q^{(t)}(r,x)= \mathbf{E}_y\left(\exp\left\{- \int_0^w K^{(t)}\left(Q^{(t)}(r-s,\xi_s^{(t)}) \right)\mathrm{d}s\right\} Q^{(t)}(r-w, \xi_{w}^{(t)}) \right),
\end{align}
where
\begin{align}
	K^{(t)}(v):= \frac{1}{v}\psi^{(t)}(v)
\end{align}
and $\psi^{(t)}(v)$ is defined in \eqref{Def-of-parameters-2}.  Also by the Markov property of $\xi^{(t)}$, we see that for all $y\in \R$ and $w\in [0,r]$, it holds that
\begin{align}
	& \Upsilon_w:= \exp\left\{- \int_0^w K^{(t)}\left(Q^{(t)}(r-s,\xi_s^{(t)}) \right)\mathrm{d}s\right\} Q^{(t)}(r-w, \xi_{w}^{(t)})\nonumber\\
	& = \mathbf{E}_y\left( \exp\left\{- \int_0^r K^{(t)}\left(Q^{(t)}(r-s,\xi_s^{(t)}) \right)\mathrm{d}s\right\} Q^{(t)}(0, \xi_{r}^{(t)}) \big| \xi_s^{(t)}: s\leq w\right).
\end{align}
Therefore, $\{(\Upsilon_w)_{w\in [0,r]}, \mathbf{P}_y\}$ is a non-negative martingale, which implies that for any stopping time $S$,
\begin{align}
		& Q^{(t)}(r,x) = \mathbf{E}_x	 (\Upsilon_{w\land S})\nonumber\\
		&= \mathbf{E}_x\left(\exp\left\{-	\int_0^{w\land S} K^{(t)}\left(Q^{(t)}(r-s,\xi_s^{(t)}) \right)\mathrm{d}s\right\}
		Q^{(t)}(r-w\land S, \xi_{w\land S}^{(t)}) \right).
\end{align}
In particular,
set $S=\tau_{-q/2}^{(t), +}:= \inf\left\{r>0: \xi_r^{(t)}\geq -q/2\right\} $
and $r=w=\delta$, we see that
\begin{align}\label{upbound1}
Q^{(t)}(\delta, -q)&\leq  \mathbf{E}_{-q}\left(Q^{(t)}(\delta-\delta\land \tau_{-q/2}^{(t),+}, \xi_{\delta \land \tau_{-q/2}^{(t),+}}^{(t)}) \right)\nonumber\\
	& = \mathbf{E}_{-q}\left(Q^{(t)}\left(\delta-\tau_{-q/2}^{(t),+}, \xi_{ \tau_{-q/2}^{(t),+}}^{(t)}\right) 1_{\{ \tau_{-q/2}^{(t),+} \leq \delta\}} \right),
\end{align}
where in the last equality we used the fact that on the event $\{\delta < \tau_{-q/2}^{(t),+}\}= \left\{\sup_{s\leq \delta} \xi_s^{(t)} <-q/2  \right\}$,  it holds that
\begin{align}
	Q^{(t)}(\delta-\delta\land \tau_{-q/2}^{(t),+}, \xi_{\delta \land \tau_{-q/2}^{(t),+}}^{(t)})= Q^{(t)}(0,  \xi_{\delta }^{(t)})= t^{\frac{1}{\alpha-1}}1_{\{\xi_{\delta }^{(t)} >0 \}}=0.
\end{align}
Note that $Q^{(t)}(\delta, -q)\leq t^{\frac{1}{\alpha-1}}$ . Note  also that, by \eqref{Tail-of-M}, for any $z<-q/4$ and $r>0$,
\[
Q^{(t)}(r, z) \leq t^{\frac{1}{\alpha-1}}\P_0\left(M>-z\sqrt{t}\right)\leq t^{\frac{1}{\alpha-1}}\P_0\left(M>q\sqrt{t}/4\right) \lesssim
q^{-\frac{2}{\alpha-1}}.
\]
Comparing $\xi_{ \tau_{-q/2}^{(t),+}}^{(t)}$ with $-q/4$, using \eqref{upbound1} and the two facts above, we get that
\begin{align}
	&Q^{(t)}(\delta, -q) \lesssim t^{\frac{1}{\alpha-1}}\mathbf{P}_{-q}\left( \xi_{ \tau_{-q/2}^{(t),+}}^{(t)}>-\frac{q}{4} \right)+ q^{-\frac{2}{\alpha-1}}
 \mathbf{P}_{-q}\left(\delta\geq \tau_{-q/2}^{(t),+} \right)\nonumber\\
	& = t^{\frac{1}{\alpha-1}}\mathbf{P}_{-q\sqrt{t}/2}\left( \xi_{ \tau_{0}^+}>\frac{q\sqrt{t}}{4} \right)+
   q^{-\frac{2}{\alpha-1}} \mathbf{P}_{0}\left(\sup_{s\leq t\delta} \xi_s\geq q\sqrt{t}/2 \right).
\end{align}
Now combining the above inequality, Lemma \ref{Moment-of-overshoot} and Doob's maximal inequality, we get that
\begin{align}
	&Q^{(t)}(\delta, -q)\lesssim \frac{t^{\frac{1}{\alpha-1}}}{(q\sqrt{t})^{r_0-2}} \sup_{x>0}\mathbf{E}_{-x}\left( \xi_{ \tau_{0}^+}^{r_0-2}\right)+ \frac{1}{q^{\frac{2}{\alpha-1}+2}t}\mathbf{E}_{0}\left( |\xi_{t\delta}|^2\right)\nonumber\\
	&\lesssim \frac{1}{q^{r_0-2} t^{\frac{r_0-2}{2}- \frac{1}{\alpha-1}}} + \frac{\delta}{q^{\frac{2\alpha}{\alpha-1}}}.
\end{align}
Since $\frac{r_0-2}{2}-\frac{1}{\alpha-1} >0$ under
{\bf (H4')},
letting $T(\delta)$ be sufficiently large so that $t^{\frac{r_0-2}{2}- \frac{1}{\alpha-1}} >\delta^{-1}$ for all $t>T(\delta)$, we get the desired \eqref{toprove}.
\hfill$\Box$
\bigskip

\textbf{Proof of Theorem \ref{thm2}: }  By Theorem \ref{thm1}, for any $\theta>0$, we have the following lower bound for the $ \liminf$:
\begin{align}\label{Lower-bound-of-local}
	 & \liminf_{t\to\infty} t^{\frac{1}{\alpha-1}}\P_{\sqrt{t}y}\left(Z_t(A)>0 \right)\geq  \lim_{t\to\infty} t^{\frac{1}{\alpha-1}}\left(1- \E_{\sqrt{t}y}\left(\exp\left\{-\frac{\theta}{t^{\frac{1}{\alpha-1}-\frac{1}{2}}}Z_t(A)\right\}\right)\right) \nonumber\\
	 &= -\log \E_{\delta_y}\left(\exp\left \{-\theta \ell(A)
 Y_1(0)\right\}\right)\stackrel{\theta\uparrow \infty}{\longrightarrow} -\log \P_{\delta_y}\left( Y_1(0)=0\right).
\end{align}
Now we prove the $ \limsup$ is no larger than the right-hand side above.
By the branching property, for
 $\kappa>0$, it holds that
\begin{align}
	\P_{\sqrt{t}y}\left(Z_{(1+\kappa)t}(A)>0 \right) = \E_{\sqrt{t}y}\left(1- \exp\left\{ \int \log \P_a(Z_{t\kappa}(A)=0) Z_t(\mathrm{d}a)\right\}\right).
\end{align}
Noticing that for all $a\in\R$,
\[
\P_a(Z_{t\kappa}(A)=0)\geq
    \P_a(Z_{t\kappa}(\R)=0)=\P_0(Z_{t\kappa}(\R)=0)
\stackrel{t\to\infty}{\longrightarrow}1.
\]
Using the fact that $\log x\sim x-1$ as $x\to 1$,
we see that  there exists $N(\kappa)$ such that as $t>N(\kappa)$,
\begin{align}\label{step_23}
	& \P_{\sqrt{t}y}\left(Z_{(1+\kappa)t}(A)>0 \right) \leq  \E_{\sqrt{t}y}\left(1- \exp\left\{ -\frac{1}{2} \int  \P_a(Z_{t\kappa}(A)>0) Z_t(\mathrm{d}a)\right\}\right).
\end{align}
 Now we fix a small $\varepsilon>0$.
Suppose that $t$ is large enough such that $A\subset [-\varepsilon\sqrt{t},\varepsilon\sqrt{t}].$  For $a< -2\varepsilon\sqrt{t}$, by Lemma \ref{Technical-Lemma}, when $t$ is large enough, we have
\begin{align}\label{step_20}
	\P_a(Z_{t\kappa}(A)>0)\leq \P_{-2\varepsilon \sqrt{t}}(M_{t\kappa}>-\varepsilon\sqrt{t})\leq \frac{C(\varepsilon)\kappa}{t^{\frac{1}{\alpha-1}}}.
\end{align}
Similarly, for $a>2\varepsilon \sqrt{t}$, when $t$ is sufficient large, it holds that
\begin{align}\label{step_21}
	\P_a(Z_{t\kappa}(A)>0)\leq
	\P_{2\varepsilon \sqrt{t}}(M_{t\kappa}^-\leq \varepsilon\sqrt{t})
	\leq \frac{C(\varepsilon)\kappa}{t^{\frac{1}{\alpha-1}}}.
\end{align}
When $|a|\leq 2\varepsilon \sqrt{t}$, by \eqref{Survival-prob-zeta}, we have
\begin{align}\label{step_22}
	\P_a(Z_{t\kappa}(A)>0) \leq \P_0\left(Z_{t\kappa}(\R)>0\right)\leq \frac{C_*}{(t\kappa)^{\frac{1}{\alpha-1}}}
\end{align}
for some constant $C_*\in (0,\infty)$. Combining \eqref{step_23}, \eqref{step_20}, \eqref{step_21}, \eqref{step_22} and the inequality $1-e^{-|x|-|y|}\leq (1-e^{-|x|}) + |y|$, we obtain that
\begin{align}\label{step_24}
	& t^{\frac{1}{\alpha-1}}\P_{\sqrt{t}y}\left(Z_{(1+\kappa)t}(A)>0 \right)\leq t^{\frac{1}{\alpha-1}} \E_{\sqrt{t}y}\left(1- \exp\left\{- \frac{C_*}{2\kappa^{\frac{1}{\alpha-1}}t^{\frac{1}{\alpha-1}} } Z_t([-2\varepsilon\sqrt{t},2\varepsilon\sqrt{t}])\right\} \right)\nonumber\\
	&\quad + \frac{C(\varepsilon)\kappa}{2}\E_{\sqrt{t}y} (Z_t((-\infty, -2\varepsilon\sqrt{t})\cup (2\varepsilon\sqrt{t},\infty )))\nonumber\\
	& \leq t^{\frac{1}{\alpha-1}} \E_{\sqrt{t}y}\left(1- \exp\left\{- \frac{C_*}{2\kappa^{\frac{1}{\alpha-1}}t^{\frac{1}{\alpha-1}} } Z_t([-2\varepsilon\sqrt{t},2\varepsilon\sqrt{t}])\right\} \right) + \frac{C(\varepsilon)\kappa}{2},
\end{align}
where in the last inequality we used the fact that $\E_{\sqrt{t}y}(Z_t(\R))=1$.  This is because, thanks to the  branching property and the criticality of our branching L\'{e}vy process, $(Z_t(\R))_{t\geq 0}$ is a martingale under $\P_x$ with $Z_0(\R)=1$ $\P_x$-a.s. for any $x\in \R$.
Now define a  function
\[
f(x):
=  \frac{C_*}{2\kappa^{\frac{1}{\alpha-1}}}\left(\left(2- \frac{1}{2\varepsilon}|x|\right)_+ \land 1\right),
x\in\R,
\]
then we see $f$ is a bounded continuous function with support equal to $[-4\varepsilon, 4\varepsilon]$
and $f=\frac{C_*}{2\kappa^{\frac{1}{\alpha-1}}}$ for $x\in [-2\varepsilon,2\varepsilon]$.
Plugging this observation into \eqref{step_24},
we get that for large $t$,
\begin{align}\label{step_25}
	&t^{\frac{1}{\alpha-1}}\P_{\sqrt{t}y}\left(Z_{(1+\kappa)t}(A)>0 \right)\nonumber\\
	& \leq t^{\frac{1}{\alpha-1}} \E_{\sqrt{t}y}\left(1- \exp\left\{- \frac{1}{t^{\frac{1}{\alpha-1}}} \int f\left(\frac{a}{\sqrt{t}}\right)Z_t(\mathrm{d} a)\right\} \right) + \frac{C(\varepsilon)\kappa}{2}\nonumber\\
	&= t^{\frac{1}{\alpha-1}} \E_{\sqrt{t}y}\left(1- \exp\left\{- \int f\left(a\right)Z_1^{(t)}(\mathrm{d} a)\right\} \right) + \frac{C(\varepsilon)\kappa}{2}.
\end{align}
Letting $t\to\infty$ in \eqref{step_25}, using Theorem \ref{thm1} with $h=0$ and $g=f$,
we get that
\begin{align}
	& \limsup_{t\to\infty} t^{\frac{1}{\alpha-1}}\P_{\sqrt{t}y}\left(Z_{t}(A)>0 \right)
  = (1+\kappa)^{\frac{1}{\alpha-1}}\limsup_{t\to\infty} t^{\frac{1}{\alpha-1}}\P_{\sqrt{t}\sqrt{1+\kappa}y}\left(Z_{(1+\kappa)t}(A)>0 \right)\nonumber\\
	& \leq (1+\kappa)^{\frac{1}{\alpha-1}} \lim_{t\to\infty}  t^{\frac{1}{\alpha-1}} \E_{\sqrt{t}\sqrt{1+\kappa}y}\left(1- \exp\left\{- \int f\left(a\right)Z_1^{(t)}(\mathrm{d} a)\right\} \right) +  (1+\kappa)^{\frac{1}{\alpha-1}}\frac{C(\varepsilon)\kappa}{2}\nonumber\\
	 & = -(1+\kappa)^{\frac{1}{\alpha-1}}   \log  \E_{\delta_{\sqrt{1+\kappa}y}}\left( \exp\left\{- \int  f(a)X_1(\mathrm{d} a)\right\} \right) + \frac{C(\varepsilon)(1+\kappa)^{\frac{1}{\alpha-1}}  }{2}\kappa
	 \nonumber\\
	 &\leq  -(1+\kappa)^{\frac{1}{\alpha-1}}   \log  \P_{\delta_{\sqrt{1+\kappa}y}}\left( X_1([-4\varepsilon, 4\varepsilon])=0\right) + \frac{C(\varepsilon)(1+\kappa)^{\frac{1}{\alpha-1}}  }{2}\kappa,
\end{align}
where in the last inequality we used the fact that $f$ is supported in $[-4\varepsilon, 4\varepsilon]$.
Letting $\kappa\to0$ in the above inequality, we see that
\begin{align}
	& \limsup_{t\to\infty} t^{\frac{1}{\alpha-1}}\P_{\sqrt{t}y}\left(Z_{t}(A)>0 \right)
	\leq
	 - \log  \P_{\delta_y}\left( X_1([-4\varepsilon, 4\varepsilon])=0\right) = - \log  \P_{\delta_y}\left( Y_1(x)=0,\ \forall |x|\leq 4\varepsilon\right) .
\end{align}
Taking $\varepsilon\to 0$ and using Remark \ref{remark}, we conclude that
\begin{align}\label{Upper-bound}
	 \limsup_{t\to\infty} t^{\frac{1}{\alpha-1}}\P_{\sqrt{t}y}\left(Z_{t}(A)>0 \right)\leq - \log  \P_{\delta_y}\left( Y_1(0)=0\right) .
\end{align}
Combining \eqref{Lower-bound-of-local} and \eqref{Upper-bound}, we complete the proof of Theorem \ref{thm2}.

\hfill$\Box$

\textbf{Proof of Theorem \ref{thm4}: }
(i) For any $f\in C_c^+(\R)$, it holds that
\begin{align}\label{step_36}
	&\E_{\sqrt{t}y}\left(\exp\left\{-\frac{1}{t^{\frac{1}{\alpha-1}-\frac{1}{2}}}\int f(x)Z_t(\mathrm{d}x)\right\}| Z_t(A)>0 \right)\nonumber\\
	& = 1- \frac{1}{\P_{\sqrt{t}y}(Z_t(A)>0)}\E_{\sqrt{t}y}\left( \left(1- \exp\left\{-\frac{1}{t^{\frac{1}{\alpha-1}-\frac{1}{2}}}\int f(x)Z_t(\mathrm{d}x)\right\}\right)1_{\{Z_t(A)>0\}}\right).
\end{align}
Let $B=(a,b)$ be a bounded interval such that $\mbox{supp}(f)\subset B$ and $A\subset B$. Then by Theorem \ref{thm2}, we see that
\begin{align}\label{step_37}
	& \left| \frac{1}{\P_{\sqrt{t}y}(Z_t(A)>0)}\E_{\sqrt{t}y}\left( \left(1- \exp\left\{-\frac{1}{t^{\frac{1}{\alpha-1}-\frac{1}{2}}}\int f(x)Z_t(\mathrm{d}x)\right\}\right)\left(1_{\{Z_t(B)>0\}} - 1_{\{Z_t(A)>0\}}\right)\right)\right| \nonumber\\
	& \leq \frac{1}{\P_{\sqrt{t}y}(Z_t(A)>0)} \left(\P_{\sqrt{t}y}(Z_t(B)>0)- \P_{\sqrt{t}y}(Z_t(A)>0)\right) \stackrel{t\to\infty}{\longrightarrow}0.
\end{align}
Further,
combining Theorem \ref{thm1} and Theorem \ref{thm2}, we get that
\begin{align}\label{step_38}
	& \lim_{t\to\infty} \frac{1}{\P_{\sqrt{t}y}(Z_t(A)>0)}\E_{\sqrt{t}y}\left( \left(1- \exp\left\{-\frac{1}{t^{\frac{1}{\alpha-1}-\frac{1}{2}}}\int f(x)Z_t(\mathrm{d}x)\right\}\right)1_{\{Z_t(B)>0\}}\right)\nonumber\\
	& = \lim_{t\to\infty} \frac{1}{\P_{\sqrt{t}y}(Z_t(A)>0)}\E_{\sqrt{t}y}\left(1- \exp\left\{-\frac{1}{t^{\frac{1}{\alpha-1}-\frac{1}{2}}}\int f(x)Z_t(\mathrm{d}x)\right\}\right)\nonumber\\
	& =   \frac{1}{\log \P_{\delta_y}(Y_1(0)=0)}\log \E_{\delta_y}\left(\exp\left\{-Y_1(0)\int f(x)\mathrm{d}x \right\}\right).
\end{align}
Therefore, combining \eqref{step_36}, \eqref{step_37} and \eqref{step_38} we see that
\begin{align}\label{step_42}
	&\lim_{t\to\infty}  \E_{\sqrt{t}y}
	\left(\exp\left\{-\frac{1}{t^{\frac{1}{\alpha-1}-\frac{1}{2}}}\int f(x)Z_t(\mathrm{d}x)\right\}| Z_t(A)>0 \right)\nonumber\\
	& = 1-  \frac{1}{\log \P_{\delta_y}(Y_1(0)=0)}\log \E_{\delta_y}\left(\exp\left\{-Y_1(0)\int f(x)\mathrm{d}x \right\}\right).
\end{align}
For any $\theta,\varepsilon>0$, taking $f= h+\frac{\theta}{2\varepsilon}1_{[-\varepsilon, \varepsilon]}$ in \eqref{N-measure-equation}  and letting $\varepsilon\to 0$, by Lemma \ref{Lemma:Density-under-N-measure}, we have
\begin{align}\label{step_43}
	&\N_y\left(1- e^{-\theta Y_1(0) - w_1(h)}\right)= \lim_{\varepsilon\to 0} \N_y\left(1- \exp\left\{-\frac{\theta}{2\varepsilon} w_1 ([-\varepsilon,\varepsilon]) -w_1(h)\right\}\right) \nonumber\\
	&= -\lim_{\varepsilon\to 0} \log \E_{\delta_y}\left(\exp\left\{-\frac{\theta}{2\varepsilon}X_1([-\varepsilon,\varepsilon]) -X_1(h)\right\}\right) \nonumber\\
	&= -\log \E_{\delta_y}\left(\exp\left\{-\theta Y_1(0) -X_1(h)\right\}\right).
\end{align}
Letting $h=0$ and $\theta\uparrow \infty$, we also see that
\begin{align}\label{step_44}
	\N_y(Y_1(0)>0) =  -\log \E_{\delta_y}\left(Y_1(0)=0\right).
\end{align}
Combining \eqref{step_42}, \eqref{step_43} and \eqref{step_44}, we have that
\begin{align}
	&\lim_{t\to\infty} \E_{\sqrt{t}y}\left(\exp\left\{-\frac{1}{t^{\frac{1}{\alpha-1}-\frac{1}{2}}}\int f(x)Z_t(\mathrm{d}x)\right\}| Z_t(A)>0 \right)\nonumber\\
	& = 1- \frac{1}{\N_y(Y_1(0)>0) }\N_y\left(1- \exp\left\{-Y_1(0)\int f(x)\mathrm{d}x\right\} \right)\nonumber\\
	& = 1- \frac{1}{\N_y(Y_1(0)>0) }\N_y\left(\left(1- \exp\left\{-Y_1(0)\int f(x)\mathrm{d}x\right\}\right)1_{\{Y_1(0)>0\}} \right)\nonumber\\
	& = \N_y\left(\exp\left\{-Y_1(0)\int f(x)\mathrm{d}x\right\} \big| Y_1(0)>0 \right),
\end{align}
which implies (i).

(ii) To prove the convergence in distribution in the weak topology, it suffices to prove that for any $g\in B_{Lip}^+(\R)$ (for example, see \cite[Lemma 3.4]{HRS}),
\begin{align}\label{Goal-thm4}
	& \lim_{t\to\infty} \E_{\sqrt{t}y}\left(\exp\left\{ - \frac{1}{t^{\frac{1}{\alpha-1}}}\int g\left(\frac{y}{\sqrt{t}}\right)Z_t(\mathrm{d}y)\right\}\big| Z_t(A)>0\right)\nonumber\\
	&= \N_y\left(\exp\left\{-w_1(g)\right\}| Y_1(0)>0\right).
\end{align}
Note that
\begin{align}\label{step_47}
	& \E_{\sqrt{t}y}\left(\exp\left\{ - \frac{1}{t^{\frac{1}{\alpha-1}}}\int g\left(\frac{y}{\sqrt{t}}\right)Z_t(\mathrm{d}y)\right\}\big| Z_t(A)>0\right) \nonumber\\
	& = \frac{1}{\P_{\sqrt{t}y}(Z_t(A)>0)}  \E_{\sqrt{t}y}\left(\exp\left\{ - \frac{1}{t^{\frac{1}{\alpha-1}}}\int g\left(\frac{y}{\sqrt{t}}\right)Z_t(\mathrm{d}y)\right\} 1_{\left\{Z_t(A)>0\right\}}\right).
\end{align}
Since $1_{\{|x|>0\}}\geq 1-e^{- a |x|}$ for $a\geq 0$, by Theorem \ref{thm1} and Theorem \ref{thm2}, for any $\theta \in (0,\infty)$, we have that
\begin{align}
	&\limsup_{t\to\infty} t^{\frac{1}{\alpha-1}}\bigg|  \E_{\sqrt{t}y}\left(\exp\left\{ - \frac{1}{t^{\frac{1}{\alpha-1}}}\int g\left(\frac{y}{\sqrt{t}}\right)Z_t(\mathrm{d}y)\right\} 1_{\left\{Z_t(A)>0\right\}}\right) \nonumber\\
	&\quad -\E_{\sqrt{t}y}\left(\exp\left\{ - \frac{1}{t^{\frac{1}{\alpha-1}}}\int g\left(\frac{y}{\sqrt{t}}\right)Z_t(\mathrm{d}y)\right\} \left(1-\exp\left\{-\frac{\theta}{t^{\frac{1}{\alpha-1}-\frac{1}{2}}}Z_t(A)\right\}\right)\right)  \bigg|\nonumber\\
	& \leq \lim_{t\to\infty} t^{\frac{1}{\alpha-1}}\bigg| \P_{\sqrt{t}y}(Z_t(A)>0)- \E_{\sqrt{t}y}\left(1- \exp\left\{-\frac{\theta}{t^{\frac{1}{\alpha-1}-\frac{1}{2}}}Z_t(A)\right\}\right) \bigg|\nonumber\\
	& =  \left|-\log \P_{\delta_y} (Y_1(0)=0)+ \log \E_{\delta_y}(\exp\left\{-\theta 	\ell(A)	Y_1(0)\right\})\right|=: G(\theta).
\end{align}
Therefore, combining the above inequality and Theorem \ref{thm2}, we conclude that for any $\theta\in (0,\infty)$,
\begin{align}\label{step_45}
	& \limsup_{t\to\infty} t^{\frac{1}{\alpha-1}}   \E_{\sqrt{t}y}\left(\exp\left\{ - \frac{1}{t^{\frac{1}{\alpha-1}}}\int g\left(\frac{y}{\sqrt{t}}\right)Z_t(\mathrm{d}y)\right\} 1_{\left\{Z_t(A)>0\right\}}\right)\nonumber\\
	&\leq G(\theta)- \log \E_{\delta_y}\left(\exp\left\{-\theta
	\ell(A)
	 Y_1(0)-X_1(g)\right\}\right) + \log \E_{\delta_y}\left(\exp\left\{-X_1(g)\right\}\right) \nonumber\\
	&\stackrel{\theta\uparrow\infty}{\longrightarrow} -\log \E_{\delta_y}\left(\exp\left\{-X_1(g)\right\}1_{\{Y_1(0)=0\}}\right) +  \log \E_{\delta_y}\left(\exp\left\{-X_1(g)\right\}\right) .
\end{align}
Similarly,  we also have that
\begin{align}\label{step_46}
	&\liminf_{t\to\infty} t^{\frac{1}{\alpha-1}}   \E_{\sqrt{t}y}\left(\exp\left\{ - \frac{1}{t^{\frac{1}{\alpha-1}}}\int g\left(\frac{y}{\sqrt{t}}\right)Z_t(\mathrm{d}y)\right\} 1_{\left\{Z_t(A)>0\right\}}\right)\nonumber\\
	&\geq -G(\theta)- \log \E_{\delta_y}\left(\exp\left\{-\theta
	\ell(A)
	 Y_1(0)-X_1(g)\right\}\right) + \log \E_{\delta_y}\left(\exp\left\{-X_1(g)\right\}\right) \nonumber\\
	&\stackrel{\theta\uparrow\infty}{\longrightarrow} -\log \E_{\delta_y}\left(\exp\left\{-X_1(g)\right\}1_{\{Y_1(0)=0\}}\right) +  \log \E_{\delta_y}\left(\exp\left\{-X_1(g)\right\}\right) .
\end{align}
Combining Theorem \ref{thm2}, \eqref{step_47} \eqref{step_45} and \eqref{step_46}, we conclude that
\begin{align}\label{step_48}
	&\lim_{t\to\infty}  \E_{\sqrt{t}y}\left(\exp\left\{ - \frac{1}{t^{\frac{1}{\alpha-1}}}\int g\left(\frac{y}{\sqrt{t}}\right)Z_t(\mathrm{d}y)\right\}\big| Z_t(A)>0\right) \nonumber\\
	& = \frac{1}{-\log \P_{\delta_y}(Y_1(0)=0)}\left( -\log \E_{\delta_y}\left(\exp\left\{-X_1(g)\right\}1_{\{Y_1(0)=0\}}\right) +  \log \E_{\delta_y}\left(\exp\left\{-X_1(g)\right\}\right) \right).\quad
\end{align}
Combining \eqref{step_43}, \eqref{step_44} and \eqref{step_48}, we get that
\begin{align}
	&\lim_{t\to\infty}  \E_{\sqrt{t}y}\left(\exp\left\{ - \frac{1}{t^{\frac{1}{\alpha-1}}}\int g\left(\frac{y}{\sqrt{t}}\right)Z_t(\mathrm{d}y)\right\}\big| Z_t(A)>0\right) \nonumber\\
	& = \frac{1}{\N_y(Y_1(0)>0)} \left(\N_y\left(1- e^{-w_1(g)}1_{\{Y_1(0)=0\}}\right) - \N_y\left(1- e^{-w_1(g)}\right)\right) \nonumber\\
	& = \frac{1}{\N_y(Y_1(0)>0)}  \N_y\left(e^{-w_1(g)}1_{\{Y_1(0)>0\}}\right) = \N_y\left(e^{-w_1(g)}|Y_1(0)>0\right),
\end{align}
which implies the desired result.

\hfill$\Box$

\section{Proof of Lemma \ref{Lemma:Density-under-N-measure}}\label{Section:Lemma:Density-under-N}

\textbf{Proof of Lemma \ref{Lemma:Density-under-N-measure}:} Suppose that $Q_t$ is the transition semigroup of the super Brownian motion, i.e.
\[
Q_t(\nu_1, \mathrm{d} \nu_2):= \P_{\nu_1}(X_t\in \mathrm{d} \nu_2).
\]
Let $Q_t^\circ$ be the restriction of $Q_t$ on $\mathcal{M}_F(\R)\setminus \{\mathbf{0}\}.$ By \cite[A.41 or (8.46)]{LZ}, for any $y\in \R$, $0< r_1<..<r_m<\infty$ and $\nu_1,...,\nu_m\in\mathcal{M}_F(\R)\setminus \{\mathbf{0}\} $, we have
\begin{align}
	\N_y\left(w_{r_1}\in \mathrm{d}\nu_1,..., w_{r_m}\in \mathrm{d}\nu_m\right)= \N_y\left(w_{r_1}\in \mathrm{d}\nu_1\right) Q_{r_2-r_1}^\circ(\nu_1,\mathrm{d}\nu_2) \cdots Q_{r_m-r_{m-1}}^\circ (\nu_{m-1}, \mathrm{d}\nu_m).
\end{align}
In particular, for any $s<t$,
\begin{align}
	& \N_y\left(w_{s}\in \mathcal{M}_F(\R)\setminus \{\mathbf{0}\},  w_t\in \mathcal{A}^c \right)= \int_{\mathcal{M}_F(\R)\setminus \{\mathbf{0}\}} \N_y(w_r\in \mathrm{d}\nu_1) Q_{t-r}^\circ(\nu_1, \mathcal{A}^c)\nonumber\\
	& = \int_{\mathcal{M}_F(\R)\setminus \{\mathbf{0}\}} \N_y(w_r\in \mathrm{d}\nu_1) \P_{\nu_1} (X_{t-r}\in \mathcal{A}^c),
\end{align}
where in the last equality we used the fact that $\mathbf{0}\notin \mathcal{A}^c$.  Since $\P_{\nu_1} (X_{t-r}\in \mathcal{A})=1$ for all $\nu_1\in \mathcal{M}_F(\R)$ and $t>r$, we obtain that
\begin{align}\label{step_41}
	& \N_y\left(w_{r}\in \mathcal{M}_F(\R)\setminus \{\mathbf{0}\},  w_t\in \mathcal{A}^c \right)
=0.
\end{align}
Moreover, $w_t\in \mathcal{A}^c$
implies
 that $w_t\neq \mathbf{0}$, therefore, it must hold that $w_r\in  \mathcal{M}_F(\R)\setminus \{\mathbf{0}\}$. Therefore, by \eqref{step_41}, we get that
\begin{align}
	& \N_y\left(w_{r}\in \mathcal{M}_F(\R)\setminus \{\mathbf{0}\},  w_t\in \mathcal{A}^c \right)=
	 \N_y\left(  w_t\in \mathcal{A}^c \right)=0,
\end{align}
which implies the desired result.

\hfill$\Box$

\bigskip
\noindent
{\bf Acknowledgements:}
We thank the referees for very helpful comments and suggestions.
We thank Zhenyao Sun for helpful discussions.
 We also thank Xinxin Chen for comments.
 Part of the research for this paper was done while the third-named author was visiting Jiangsu Normal University,
where he was partially supported by a grant from the National Natural Science Foundation of China
(12471139, Yingchao Xie).
\bigskip
\noindent

\begin{singlespace}

\end{singlespace}

\vskip 0.2truein
\vskip 0.2truein

\noindent{\bf Haojie Hou:}  School of Mathematics and Statistics, Beijing Institute of Technology,   Beijing 100081, P.R. China. Email: {\texttt houhaojie@bit.edu.cn}

\smallskip

\noindent{\bf Yan-Xia Ren:} LMAM School of Mathematical Sciences \& Center for
Statistical Science, Peking
University,  Beijing, 100871, P.R. China. Email: {\texttt
yxren@math.pku.edu.cn}

\smallskip
\noindent {\bf Renming Song:} Department of Mathematics,
University of Illinois Urbana-Champagn,
Urbana, IL 61801, U.S.A.
Email: {\texttt rsong@illinois.edu}


\small
\begin{thebibliography}{99}



\bibitem{BMS2024} Barnes, C., Mytnik, L. and Sun, Z.: On the coming down from infinity of coalescing Brownian motions. \emph{Ann. Probab.} \textbf{52}(2024) 67--92.


\bibitem{BM2018} Bertoin, J. and Mallein, B.: Biggins'  martingale convergence for branching L\'{e}vy processes. \emph{Electron. Commun. Probab.} \textbf{23} (2018) 1--12.


\bibitem{BM} Bertoin, J. and Mallein, B.: Infinitely ramified point measures and branching L\'{e}vy processes. \emph{Ann. Probab.} \textbf{47} (2019) 1619--1652.


\bibitem{BCG1997}  Bramson, M. Cox, J.T. and Greven, A.: Invariant measures of critical spatial branching process in high dimensions. \emph{Ann. Probab.} 	\textbf{25}(1997) 56--70.


\bibitem{Chen}  Bai, T.,  Chen,  X. and Lin, S. :  Critical branching random walks on $\mathbb{Z}^2$: local survival probabilities and Yaglom limit theorems. arXiv: 2609.15538.


\bibitem{Durrett1979} Durrett, R.: An infinite particle system with additive interactions.
\emph{Adv. in Appl. Probab.}
\textbf{11}(1979) 355--383.

\bibitem {E.B1.}
Dynkin, E. B.: Branching exit Markov systems and superprocesses. \emph{Ann. Probab.} \textbf{29} (2001), 1833--1858.

\bibitem{DyKu} Dynkin, E. B. and Kuznetsov, S. E.: $\mathbb{N}$-measures for branching exit Markov systems and their applications to differential equations.
\emph{Probab. Theory Relat. Fields} \textbf{130}(1) (2004), 135--150.



\bibitem{Feller} Feller, W.: \emph{An Introduction to Probability Theory and Its Applications}, Vol. 2. Wiley, New York, 1964.


\bibitem{FMW2010} Fleischmann, K., Mytnik, L. and Wachtel, V.:
Optimal local H\"{o}lder index for density states of superprocesses with $(1+\beta)$-branching mechanism.
\emph{Ann. Probab.} \textbf{38}(2010) 1180--1220.


\bibitem{GX} Grama, I. and Xiao, H.: Conditioned local limit theorems for random walks on the real line.
\emph{Ann. Inst. Henri  Poincar\'{e} Probab. Stat.}
\textbf{61} (2025), 403--456.


\bibitem{HL2024} Hong, W. and Liang, S.:
Conditional central limit theorem for critical branching random walk.
\emph{ALEA Lat. Am. J. Probab. Math. Stat.} \textbf{21}(2024)  555--574.

\bibitem{HJRS} Hou, H., Jiang, Y., Ren, Y.-X. and Song, R.:  Tail probability of maximal displacement in critical branching
L\'evy
process with stable branching.
\emph{Bernoulli} \textbf{31} (2025), 630--648.



\bibitem{HRS} Hou, H., Ren, Y.-X. and Song, R.: Tails of extinction time and maximal displacement for critical branching killed L\'evy process.  \emph{Potential Anal.} \textbf{63} (2025), 1811--1867.



\bibitem{KV2024} Kagan, A. and V\'{e}chambre, G.: The Leftmost particle of branching subordinators. arXiv: 2409.16617.


\bibitem{Kallenberg} Kallenberg, O.:
\emph{Foundations of modern probability.}
2nd edition, Springer-Verlag, New York, 2002. xx+638 pp.


\bibitem{Kolmogorov38} Kolmogorov, A.: Zur l\"{o}sung einer biologischen aufgabe. \emph{Comm. Math. Mech. Chebyshev Univ. Tomsk.} \textbf{2}(1) (1938), 1--12.

\bibitem{KS1988} Konno, N. and Shiga, T.: Stochastic partial differential equations for some
measure-valued diffusions.
\emph{Probab. Theory Relat. Fields}
\textbf{79}(1988) 201--225.

\bibitem{LS15} Lalley, S. P.  and Shao, Y.: On the maximal displacement of critical branching random walk. \emph{Probab. Theory Relat. Fields} \textbf{162}(1--2) (2015) 71--96.

\bibitem{LS16} Lalley, S. P.  and Shao, Y.:  Maximal displacement of critical branching symmetric stable processes. \emph{Ann. Inst. Henri Poincar\'{e} Probab. Stat. } \textbf{52} (3) (2016), 1161--1177.


\bibitem{LZ11} Lalley, S. P. and Zheng, X.: Occupation statistics of critical branching random walks in two or higher dimensions. \emph{Ann. Probab.} \textbf{39}(2011) 327--368.

\bibitem{LeGall1996} Le Gall, J.-F.:  A probabilistic approach to the trace at the boundary for solutions of a semi-linear parabolic partial differential equation.  \emph{J. Appl. Math. Stoch. Anal.} \textbf{9}(1996) 399--414.

\bibitem{Lee}  Lee, T.-Y.: Conditional limit distributions of critical branching Brownian motions.  \emph{Ann. Probab.} \textbf{19} (1) (1991) 289--311.

\bibitem{LZ} Li, Z.:
\emph{Measure valued branching Markov processes.}
2nd edition, \emph{ Probab. Theory Stoch. Model.} \textbf{103} Springer, Berlin,  2022.


\bibitem{MS2023} Mallein, B. and Shi, Q.: A necessary and sufficient condition for the convergence of the derivative martingale in a branching L\'{e}vy process. \emph{Bernoulli} \textbf{29} (2023) 597--624.


\bibitem{MV1999} Marcus, M. and V\'eron, L.: Initial trace of positive solutions of some nonlinear parabolic equations. \emph{Comm. Partial Differential Equations.} \textbf{24}(1999) 1445--1499.


\bibitem{Neuman-Zheng} Neuman, E. and Zheng, X.: On the maximal displacement of near-critical branching random walks. \emph{Probab. Theory Relat. Fields} \textbf{180} (2021), 199--232.

\bibitem{Petrov1975} Petrov, V. V.:
\emph{Sums of independent random variables.}
Springer-Verlag, New York-Heidelberg, 1975. x+346 pp.

\bibitem{EP19} Powell, E.: An invariance principle for branching diffusions in bounded domains. \emph{Probab. Theory Relat. Fields} \textbf{173}(3--4) (2019), 999--1062.

\bibitem{Profeta21} Profeta, C.: Extreme values of critical and subcritical branching stable processes with positive jumps.  \emph{ALEA, Lat. Am. J. Probab. Math. Stat.} \textbf{19}(2) (2022) 1421--1433.

\bibitem{Profeta22}  Profeta, C.: Maximal displacement of spectrally negative branching L\'{e}vy processes. \emph{Bernoulli} \textbf{30} (2024) 961--982.

\bibitem{Rapenne} Rapenne, V.:
Invariant measures of critical branching random walks in high dimension.
\emph{Electron. J. Probab.} \textbf{28}(2023) 1--38.


\bibitem{RSZ2021} Ren, Y.-X., Song, R. and Zhang, R.: The extremal process of super-Brownian motion. \emph{Stoch. Proc. Appl.} \textbf{137}(2021) 1--34.


\bibitem{SW2019} Shi, Q. and  Watson, A. R.: Probability tilting of compensated fragmentations. \emph{Electronic. J. Probab.} \textbf{24}(2019) 1--39.


\bibitem{Slack1968}  Slack, R.: A branching process with mean one and possibly infinite variance. \emph{Z. Wahrscheinlichkeitstheor. Verwandte Geb.} \textbf{9} (1968) 139--145.

\bibitem{Zolotarev1957}
Zolotarev, V.: More exact statements of several theorems in the theory of branching processes. \emph{Teor. Veroyatn. Primen.} \textbf{2} (1957) 256--266.


\end{thebibliography}
\end{document}